\documentclass[a4paper, 10pt]{article}
\usepackage{geometry}
\usepackage[francais, english]{babel}
\usepackage{latexsym,amsfonts, amssymb}
\usepackage{amsthm}
\usepackage[latin1]{inputenc}
\usepackage{amsmath}
\usepackage{amssymb}
\usepackage{amsfonts}
\usepackage{dsfont}
\usepackage{hyperref}
\usepackage{amsthm}
\usepackage{latexsym}
\usepackage{epsfig}
\usepackage{mathrsfs}

\makeatletter
\renewcommand\theequation{\thesection.\arabic{equation}}
\@addtoreset{equation}{section}
\makeatother

\def\lg{\textquotedblleft}
\def\rg{\textquotedblright}

\newfont{\bbf}{cmbx12 scaled 1435}

\newtheorem{thm}{Theorem}[section]
\newtheorem{prop}{proposition}[section]
\newtheorem{lem}{Lemma}[section]

\renewcommand{\theequation}{\thesection.\arabic{equation}}
\newcommand{\eop}{\hspace*{\fill} \ensuremath{\Box}}

\begin{document}

\newcommand{\I}[1]{\mathds{1}_{{#1}}}

\def\Sum{ \displaystyle \sum }
\def\Frac{\displaystyle \frac}

\def\Cit{\mathbb{C}}
\def\esp{\mathbb{E}}
\def\Var{\hbox{\rm Var}}
\def\Cov{\hbox{\rm Cov}}
\def\Supp{\hbox{\rm Supp}}
\def\Card{\hbox{\rm Card}}
\def\Det{\hbox{\rm Det}}
\def\Tr{\hbox{\rm Tr}}
\def\Dim{\hbox{\rm dim}}
\def\Rank{\hbox{\rm dim}}
\def\Id{\hbox{\rm Id}}
\def\Ker{\hbox{\rm Ker}}
\def\ind{\mathbb{I}}
\def\Nit{\mathbb{N}}
\def\Rit{\mathbb{R}}
\def\Zit{\mathbb{Z}}
\def\prob{\mathbb{P}}
\def\where{\rm where}
\def\with{\rm with}
\def\and{\rm and}
\def\I{\rm I}
\def\J{\rm J}
\def\Ip{\rm I^{\prime}}
\def\Jp{\rm J^{\prime}}
\def\as{\rm a.s}
\def\and{\rm and}

\begin{center}
\setcounter{page}{1} \setlength{\baselineskip}{.32in} {\Large\bf
Nonparametric kernel estimation of the probability density
function of  regression errors using estimated residuals}


\vspace*{0.25cm}
 By Rawane SAMB
\\
\it{Université Pierre et Marie Curie, LSTA. }


\end{center}
\begin{center}
{\bf Abstract}
\end{center}
\indent This paper deals with the nonparametric density estimation
of the regression error term assuming its independence with the
covariate. The difference between the feasible estimator which
uses the estimated residuals and the unfeasible one using the true
 residuals
is studied. An optimal choice of the bandwidth used to estimate
the residuals is  given. We also study the asymptotic normality of
the feasible kernel estimator and its rate-optimality.
\begin{small}
\vspace{0.17cm}
\end{small}
\renewcommand{\thefootnote}{\arabic{footnote}} \setcounter{footnote}{1}
\setlength{\baselineskip}{.26in}
\section{Introduction}

 Consider a sample
$(X,Y),(X_{1},Y_{1}),\ldots,(X_{n},Y_{n})$ of independent and
identically distributed (i.i.d) random variables, where Y is the
univariate dependent variable and the covariate X is of dimension
$d$. Let $m(\cdot)$ be the conditional expectation of $Y$ given
$X$ and let
 $\varepsilon$ be
the related regression error term, so that the regression error
model is
\begin{eqnarray}
Y_i = m(X_i)+\varepsilon_i, \quad i=1,\ldots,n. \label{Rm}
\end{eqnarray}
 We wish to estimate the probability distribution function
(p.d.f) of the regression error term, $f(\cdot)$, using the
nonparametric
 residuals.
Our potential applications are as follows.
 First, an estimation of the p.d.f of $\varepsilon$ is an
 important tool for understanding the residuals behavior and
 therefore the fit of the regression model (\ref{Rm}).
This estimation of $f(\cdot)$ can be used for
 goodness-of-fit tests of a specified error  distribution in a
 parametric regression setting. Some examples can be founded in
  Loynes (1980), Akritas and Van  Keilegom (2001),  Cheng and Sun
  (2008).
   The estimation of the density of the regression error term can
   also
 be useful for testing the symmetry of the residuals.
  See Ahmad et Li (1997), Dette et
 {\it al.} (2002). Another interest of the estimation of $f$
 is that it can be used for constructing nonparametric estimators
 for the density and hazard function of $Y$ given $X$, as related
 in Van Keilegom and Veraverbeke (2002). This estimation of $f$ is also
 important when are interested in the estimation
 of the p.d.f of the response variable $Y$. See Escanciano and
 Jacho-Chavez (2010).
  Note also that an estimation of the p.d.f of the
 regression
 errors  can be useful for proposing a mode forecast of
 $Y$ given $X=x$. This mode forecast is based on an estimation of
 $m(x)+\arg\min_{\epsilon\in\Rit} f(\epsilon)$.

\vskip 0.3cm
 Relatively little is known about the nonparametric
estimation of the p.d.f and  the cumulative distribution function
(c.d.f) of the regression error. Up to few
 exceptions, the nonparametric literature focuses on studying the
 distribution of $Y$ given $X$. See Roussas
(1967, 1991), Youndj\'e (1996) and
 references therein.
Akritas and Van Keilegom (2001)  estimate the cumulative
distribution function of the regression error in heteroscedastic
model. The estimator proposed by these authors is based on a
nonparametric estimation of the residuals. Their result show the
impact of the estimation of the residuals on the limit
distribution of the
 underlying estimator
of the cumulative distribution function. M\"{u}ller, Schick and
Wefelmeyer (2004) consider the estimation of
 moments of the regression error.
Quite surprisingly, under appropriate conditions, the estimator
based
 on the true errors
is less efficient than the  estimator which uses the nonparametric
estimated residuals.
 The reason is that the latter estimator  better uses the fact that
 the regression error $\varepsilon$
has mean zero. Efromovich (2005) consider adaptive estimation of
the p.d.f of the
 regression error.
 He gives a nonparametric estimator based on the
estimated residuals, for which the Mean Integrated Squared Error
(MISE)  attains the  minimax rate. Fu and Yang (2008) study the
asymptotic normality of the estimators of the regression error
p.d.f  in nonlinear autoregressive models.
 Cheng (2005)  establishes
the asymptotic normality of an estimator of
 $f(\cdot)$ based on the
estimated residuals. This estimator is  constructed by splitting
the
 sample into two parts: the first part is used for the construction of
 estimator of $f(\cdot)$,
 while the second part of the sample is used for the estimation of the
 residuals.

\vskip 0.3cm
 The focus of this paper is to estimate  the p.d.f
of the regression error using the estimated residuals,  under the
assumption that the covariate $X$ and the regression error
 $\varepsilon$ are independent.
 In a such setup, it would be unwise to use a conditional approach
 based
on the fact that
$f(\epsilon)=f(\epsilon|x)=\varphi\left(m(x)+\epsilon| x\right)$,
where $\varphi(\cdot| x)$  is the p.d.f of $Y$ given $X=x$.
Indeed,
 the estimation of $m(\cdot)$ and $\varphi(\cdot| x)$ are affected by
the curse of dimensionality, so that the resulting estimator of
$f(\cdot)$ would have considerably a slow rate of convergence if
the
 dimension of $X$ is
high. The  approach proposed here uses a two-steps procedure
which, in
 a
first step, replaces the unobserved regression error terms by some
nonparametric estimator $\widehat{\varepsilon}_i$. In a second
step,
 the
estimated $\widehat{\varepsilon}_i$'s are used to estimate
nonparametrically $f(\cdot)$, as if they were the true
 $\varepsilon_i$'s. If
proceeding so can circumvent the curse of dimensionality, a
challenging
 issue is to evaluate
the impact of the estimated  residuals on the final estimator of
$f(\cdot)$. Hence one of the contributions of our study is to
analyze the effect of the estimation of the residuals on the
regression errors p.d.f. Kernel estimators.
 Next, an optimal
choice of
 the bandwidth used to
 estimate the residuals is  given. Finally, we study the asymptotic
 normality of the
feasible Kernel estimator and its rate-optimality.

\vskip 0.3cm
 The rest of this paper is organized as follows.
Section 2 presents ours estimators
 and proposes an  asymptotic normality of the (naive) conditional
  estimator of the
 density of the regression error term. Sections 3 and 4 group our
assumptions and  main results. The conclusion of this chapter is
given in Section 5, while the proofs of our results are gathered
in section 6 and in an appendix.

\setcounter{subsection}{0} \setcounter{equation}{0}
\renewcommand{\theequation}{\thesection.\arabic{equation}}

\section{Some nonparametric conditional estimator of the density of the regression error}

 To illustrate the potential impact of the dimension $d$ of the
$X_i$'s,
 let us first consider a naive conditional estimator of the p.d.f
 $f(\cdot)$ of the regression error term $\varepsilon$. Let $\varphi
 (\cdot|x)$ and $f(\cdot|x)$ be respectively the p.d.f. of $Y$ and $\varepsilon$
 given $X=x$. Since $f(\epsilon|x) = \varphi (m(x)+\epsilon|x)$, using
 the independence of $X$ and $\varepsilon$ gives
\begin{equation}
f(\epsilon)
 =
 f(\epsilon|x)
 =
 \varphi\left(m(x)+\epsilon|x\right).
\label{Naive}
\end{equation}
Consider some Kernel functions $K_0$, $K_1$ and some bandwidths
$b_0$, $h_0$
 and $h_1$. The expression (\ref{Naive}) of $f$ suggests to use the
 Kernel nonparametric estimator
\begin{eqnarray*}
\widetilde{f}_n(\epsilon|x)
 =
\frac{ \frac{1}{nh_0^dh_1} \sum_{i=1}^n
 K_0\left(\frac{X_i-x}{h_0}\right)
 K_1\left(\frac{Y_i-\widehat{m}_n(x)-\epsilon}{h_1}\right)
 }{\frac{1}{nh_0^d}
 \sum_{i=1}^n
 K_0\left(\frac{X_i-x}{h_0}\right)}\;,
\end{eqnarray*}
where $\widehat{m}_n(x)$ is the Nadaraya-Watson (1964) estimator
of $m(x)$ defined as
\begin{eqnarray}
\widehat{m}_n(x)
 =
 \frac{
\sum_{j=1}^n Y_{j}
 K_0\left(\frac{X_j-x}{b_0}\right)}
 {
 \sum_{j=1}^n
 K_0\left(\frac{X_j-x}{b_0}\right)}\;.
 \label{mchap}
\end{eqnarray}
The first result presented in this chapter is the following
proposition.
\begin{prop}
Define
$$
\mu_1(x,\epsilon) =
 \frac{\partial^2\varphi\left(x,
m(x)+\epsilon\right)}
 {\partial^2 x}
\int z K_0(z) z^{\top} dz, \quad \mu_2(x,\epsilon) =
\frac{\partial^2\varphi\left(x, m(x)+\epsilon\right)}
 {\partial^2 y}
\int v^2 K_1(v)dv,
$$
and suppose that  $h_0$ decrease to $0$ such that
 $nh_0^{2d}/\ln n\rightarrow\infty$,
 $\ln(1/h_0)/\ln(\ln n)\rightarrow\infty$ and
$$
{(\rm\bf{A}_0):} \quad nh_0^dh_1\rightarrow\infty, \quad
\left(\frac{nh_0^d}{h_1}\right) \left(b_0^4+\frac{\ln
n}{nb_0^d}\right) =o(1),
$$
 when $n\rightarrow\infty$. Then under Assumptions  $(A_1)-(A_{10})$
  given in the next section, we have
\begin{eqnarray*}
\sqrt{nh_0^dh_1}
 \left(
\widetilde{f}_n(\epsilon|x)
 -
 \overline{\widetilde{f}}_n(\epsilon|x)
 \right)
\stackrel{d}{\rightarrow} \mathcal{N} \left( 0,
\frac{f(\epsilon|x)}{g(x)}
 \int\int
 K_0^2(z)
 K_1^2(v)
 dz dv
\right),
\end{eqnarray*}
where $g(\cdot)$ is the marginal density of $X$ and
\begin{eqnarray*}
\overline{\widetilde{f}}_n(\epsilon|x)
 =
 f(\epsilon| x)
 +
 \frac{h_0^2\mu_1(x,\epsilon)}{2g(x)}
 +
 \frac{h_1^2\mu_2(x,\epsilon)}{2g(x)}
 +
 o\left(h_0^2+h_1^2\right).
\end{eqnarray*}
\label{Prop}
\end{prop}
\noindent
 This results suggests that an optimal choice of the
 bandwidths $h_0$ and $h_1$ should achieve the minimum of
the  asymptotic mean square expansion first order terms
$$
AMSE \left( \widetilde{f}_n(\epsilon| x) \right)
 =
 \left[
 \frac{h_0^2\mu_1(x,\epsilon)}{2g(x)}
 +
 \frac{h_1^2\mu_2(x,\epsilon)}{2g(x)}
\right]^2 + \frac{ f(\epsilon|x) \int\! K_0^2(z)dz \int\!
K_1^2(v)dv} {n h_0^d h_1g(x)}\;.
$$
Elementary calculations yield that the resulting optimal
bandwidths $h_0$ and $h_1$
 are all proportional to
$n^{-1/(d+5)}$, leading to the exact consistency rate
$n^{-2/(d+5)}$ for $\widetilde{f}_n(x|\epsilon)$. In the case
$d=1$, this rate is $n^{-1/3}$, which is worst than the
 rate $n^{-2/5}$ achieved  by the optimal Kernel estimator of an univariate
 density. See  Bosq and Lecoutre (1987),
Scott (1992), Wand and Jones (1995). Note also that the exponent
$2/(d+5)$ decreases to $0$ with
 the dimension $d$. This indicates a negative impact of the dimension $d$
 on the performance of the estimator, the so-called curse of
 dimensionality.
The fact that $\widetilde{f}_n(\epsilon|x)$ is affected by the
curse
 of dimensionality is a consequence of conditioning. Indeed, (\ref{Naive})
 identifies the unconditional $f(\epsilon)$ with the conditional
 distribution
 of the regression error given the covariate.

\vskip 0.3cm To avoid this curse of dimensionality in the
nonparametric kernel estimation of $f(\epsilon)$, our approach
proposed here builds, in
 a first step,  the estimated residuals
\begin{equation}
\widehat{\varepsilon}_i = Y_i-\widehat{m}_{in}, \quad
i=1,\ldots,n, \label{epschap}
\end{equation}
where $\widehat{m}_{in} =\widehat{m}_{in}(X_i)$ is a leave-one out
 version of the Kernel regression estimator (\ref{mchap}),
\begin{equation}
\widehat{m}_{in} = \frac{\sum_{j=1\atop j\neq i}^nY_j
K_0\left(\frac{X_j-X_i}{b_0}\right)} {\sum_{j=1\atop j\neq i}^n
K_0\left(\frac{X_j-X_i}{b_0}\right)}\;. \label{mchapi}
\end{equation}
It is tempting to use, in a second step, the estimated
 $\widehat{\varepsilon}_i$ as if they were the true residuals $\varepsilon_i$. This
 would ignore that the $\widehat{m}_{n}(X_i)$'s can deliver  severely
 biased estimations of the $m(X_i)$'s for those $X_i$ which are close to the
 boundaries of the support $\mathcal{X}$ of the covariate distribution.
 To that aim, our proposed estimator trims the observations $X_i$
 outside an inner subset $\mathcal{X}_0$ of
$\mathcal{X}$,
\begin{equation}
\widehat{f}_{1n}(\epsilon)
 =
\frac{1} {b_1\sum_{i=1}^n \mathds{1}
\left(X_i\in\mathcal{X}_0\right)} \sum_{i=1}^n \mathds{1} \left(
X_i\in \mathcal{X}_0\right)
K_1\left(\frac{\widehat{\varepsilon}_i-\epsilon}{b_1}\right).
 \label{fnchap}
\end{equation}
This estimator is the so-called two-steps Kernel estimator of
$f(\epsilon)$. In principle, it would be possible to assume that
$\mathcal{X}_0$ grows
 to $\mathcal{X}$ with a negligible rate compared to the bandwidth
 $b_1$.
This would give an estimator close to the more natural Kernel
estimator $\sum_{i=1}^n
 K\left((\widehat{\varepsilon}_i-\epsilon)/b_1\right)/(nb_1)$.
However, in the rest of the paper, a fixed subset $\mathcal{X}_0$
will
 be considered for the sake of simplicity.

Observe that the two steps Kernel estimator
$\widehat{f}_{1n}(\epsilon)$ is a feasible estimator in the sense
that it does not depend on any unknown
 quantity, as desirable in practice. This contrasts with the unfeasible
 ideal Kernel estimator
\begin{equation}
 \widetilde{f}_{1n}(\epsilon)
 =
\frac{1} {b_1\sum_{i=1}^n \mathds{1}
\left(X_i\in\mathcal{X}_0\right)} \sum_{i=1}^n \mathds{1}\left(
X_i \in \mathcal{X}_0\right)
 K_1\left(\frac{\varepsilon_i-\epsilon}{b_1}\right),
 \label{fn}
\end{equation}
 which depends in particular on the unknown regression error
terms. It
 is however intuitively clear that $\widehat{f}_{1n}(\epsilon)$ and
 $\widetilde{f}_{1n}(\epsilon)$ should be closed, as illustrated by the results of the
next section.

\section{Assumptions}

The following assumptions are used in our mains results.
 \vskip 0.3cm \noindent
 {$\bf(A_1)$}
 {\it
 The support $\mathcal{X}$ of $X$  is a  compact subset of
 $\Rit^d$
and $\mathcal{X}_0$
 is an inner closed subset of $\mathcal{X}$ with non empty
 interior,
}
\medskip
\\{$\bf(A_2)$}
{\it the p.d.f. $g(\cdot)$ of the i.i.d. covariates $X, X_i$  is
strictly positive over $\mathcal{X}_0$, and has continuous second
order partial derivatives  over $\mathcal{X}$, }
\medskip
\\{$\bf(A_3)$}
{\it
 the regression function $m(\cdot)$ has continuous second order partial
derivatives  over  $\mathcal{X}$, }
\medskip
\\{$\bf(A_4)$}
{\it the i.i.d. centered error regression terms $\varepsilon,
 \varepsilon_i$'s, have finite 6th moments, and are independent of the
 covariates $X,X_i$'s,
}
\\{$\bf(A_5)$}
{\it
 the  probability  density function $f(\cdot)$ has  bounded continuous
 second order
 derivatives over  $\Rit$ and satisfies, for $h_p (e) = e^p f(e)$,
 $\sup_{e\in\Rit}|h_p^{(k)}(e)|<\infty$,
 $p\in[0,2]$, $k\in[0,2]$,
  }
\medskip
\\{$\bf(A_6)$}
{\it the p.d.f  $\varphi$ of $(X, Y)$  has bounded continuous
second  order partial derivatives over $\Rit^d\times\Rit$, }
\\{$\bf(A_7)$}
{\it the Kernel function $K_0$  is  symmetric, continuous over
$\Rit^d$ with support  contained in $[-1/2, 1/2]^d$ and $\int\!
K_0 (z) dz = 1$, }
\medskip
\\{$\bf(A_8)$}
{\it
 the Kernel function $K_{1}$ has a compact support,  is three times
 continuously differentiable over
 $\Rit$, and satisfies $\int\! K_1 (v) dv = 1$ and
$\int\! v K_1 (v) dv = 0$, }
\medskip
\\{$\bf(A_{9})$}
{\it
 the bandwidth $b_0$ decreases to $0$ and satisfies,
 for $d^*=\sup\{d+2,2d\}$, $nb_0^{d^*}/\ln n\rightarrow\infty$ and
 $\ln(1/b_0)/\ln(\ln n)\rightarrow\infty$
 when $n\rightarrow\infty$,
}
\\
{$\bf(A_{10})$}
 {\it
the bandwidth $b_1$ decreases to $0$ and satisfies
$n^{(d+8)}b_1^{7(d+4)}\rightarrow\infty$ when
$n\rightarrow\infty$. }

\vskip 0.3cm \noindent
 Assumptions $(A_2)$, $(A_3)$, $(A_5)$ and $(A_6)$ impose that all the
functions
 to be estimated nonparametrically have two bounded derivatives.
 Consequently the conditions  $\int\! z K_0 (z) dz = 0$ and
$\int\!v K_1 (v) dv = 0$, as assumed in $(A_7)$ and $(A_8)$,
represent
 standard conditions ensuring  that the bias of the resulting
 nonparametric estimators (\ref{mchap})
and (\ref{fn}) are of order $b_0^2$ and $b_1^2$.
 Assumption $(A_4)$ states independence between
the regression error terms
 and the covariates, which is the main condition for (\ref{Naive}) to
 hold.
The differentiability of $K_1$ imposed in $(A_8)$ is more specific
to our
 two-steps estimation method. Assumption $(A_8)$ is used to expand the
 two-steps Kernel estimator
$\widehat{f}_{1n}$ in (\ref{fnchap}) around the unfeasible one
$\widetilde{f}_{1n}$ from
 (\ref{fn}), using the residual error estimation
 $\widehat{\varepsilon}_i - \varepsilon_i$'s and the derivatives of $K_1$
 up to third order.
Assumption $(A_9)$ is useful for obtaining the uniform convergence
 of the regression estimator $\widehat{m}_n$ defined in (\ref{mchap}) (see for
 instance Einmahl and Mason, 2005), and also gives a similar consistency
 result for the leave-one-out estimator $\widehat{m}_{in}$ in
 (\ref{mchapi}). Assumption $(A_{10})$ is needed in the study of
  the difference between the
 feasible estimator $\widehat{f}_{1n}$ and the unfeasible
  estimator $\widetilde{f}_{1n}$.

\section{Main results}
This section is devoted to our main results. The first result we
give here concerns the pointwise consistency of the nonparamatric
Kernel estimator $\widehat{f}_{1n}$ of the density $f$. Next, the
optimal first-step and second-step bandwidths used to estimated
$f$ are proposed. We finish this section by establishing an
asymptotic normality for the estimator $\widehat{f}_{1n}$.

\subsection{Pointwise weak consistency}

The next  result gives the order of the difference between the
feasible estimator and the theoretical density of the regression
error at a fixed point $\epsilon$.

\begin{thm}
Under $(A_1)-(A_5)$ and $(A_7)-(A_{10})$, we have, when $b_0$ and
$b_1$ go to $0$,
$$
\widehat{f}_{1n}(\epsilon)-f(\epsilon)
 =
 O_{\prob}
\biggl(AMSE(b_1)+R_n(b_0, b_1)\biggr)^{1/2},
$$
where
$$
AMSE(b_1)
 =
 \esp_n
 \left[
\left(\widetilde{f}_{1n}(\epsilon)-f(\epsilon)\right)^2
 \right]
 =
 O_{\prob}
\left(b_1^4 + \frac{1}{nb_1}\right),
$$
and
\begin{eqnarray*}
 R_n(b_0, b_1)
 =
 b_0^4
 +
 \left[
 \frac{1}{(nb_1^5)^{1/2}}
 +
 \left(\frac{b_0^d}{b_1^3}\right)^{1/2}
 \right]^2
 \left(
 b_0^4
 +
 \frac{1}{nb_0^d}
 \right)^2
 +
 \left[
 \frac{1}{b_1}
 +
 \left(\frac{b_0^d}{b_1^7}\right)^{1/2}
 \right]^2
 \left(b_0^4+\frac{1}{nb_0^d}\right)^3.
\end{eqnarray*}
 \label{thm1}
\end{thm}

\noindent The result of Theorem \ref{thm1} is based on the
evaluation of the difference between  $\widehat{f}_{1n}(\epsilon)$
and $\widetilde{f}_{1n}(\epsilon)$.  This evaluation  gives an
indication about the impact of the estimation of the residuals on
the nonparametric estimation of the regression error density.

\subsection{Optimal first-step and second-step bandwidths
for the pointwise weak consistency}

 As shown in the next result, Theorem \ref{Optimbandw1}
gives some guidelines for the choice of the optimal bandwidth
$b_0$ used in the nonparametric regression errors estimation. As
far as we know, the choice of an optimal $b_0$ has not been
addressed before. In what follows, $a_n \asymp b_n$ means that $
a_n= O (b_n)$ and $b_n = O(a_n)$, i.e. that there is a constant
$C>0$ such that $|a_n|/C \leq |b_n| \leq C |a_n|$ for $n$ large
enough.

\begin{thm}
 Suppose that
$(A_1)-(A_5)$ and $(A_7)-(A_{10})$ are satisfied, and  define
$$
b_0^* = b_0^*(b_1)
 =
\arg\min_{b_0}
 R_n (b_0, b_1).
$$
where the minimization is performed over bandwidth $b_0$
fulfilling $(A_9)$. Then the bandwidth $b_0^*$ satisfies
$$
b_0^* \asymp \max \left\lbrace
\left(\frac{1}{n^2b_1^3}\right)^{\frac{1}{d+4}} ,
\left(\frac{1}{n^3b_1^7}\right)^{\frac{1}{2d+4}} \right\rbrace,
$$
and we have
$$
R_n(b_0^*, b_1) \asymp \max \left\lbrace
\left(\frac{1}{n^2b_1^3}\right)^{\frac{4}{d+4}} ,
\left(\frac{1}{n^3b_1^7}\right)^{\frac{4}{2d+4}}
 \right\rbrace.
$$
\label{Optimbandw1}
\end{thm}
Our next theorem gives the conditions for which the estimator
$\widehat{f}_{1n}(\epsilon)$ reaches the optimal rate $n^{-2/5}$
when $b_0$ takes the value $b_0^*$.
 We prove that for
$d\leq 2$, the bandwidth that minimizes the term
$AMSE(b_1)+R_n(b_0^*, b_1)$ has the same order as $n^{-1/5}$,
yielding the optimal order $n^{-2/5}$ for
$\left(AMSE(b_1)+R_n(b_0^*, b_1)\right)^{1/2}$.

\begin{thm}
Assume that $(A_1)-(A_5)$ and $(A_7)-(A_{10})$ are satisfied, and
set
$$
b_1^*
 =
\arg\min_{b_1} \biggl(AMSE(b_1)+R_n(b_0^*,b_1)\biggr),
$$
where $b_0^*=b_0^*(b_1)$ is defined as in Theorem
\ref{Optimbandw1}.
 Then
\begin{enumerate}
\item For $d\leq 2$, the bandwidth $b_1^*$ satisfies
$$
b_1^* \asymp \left(\frac{1}{n}\right)^{\frac{1}{5}},
$$
and we have
$$
\biggl( AMSE(b_1^*) +
 R_n(b_0^*, b_1^*)
\biggr)^{\frac{1}{2}} \asymp
\left(\frac{1}{n}\right)^{\frac{2}{5}}.
$$
\item For $d\geq 3$, $b_1^*$ satisfies
$$
b_1^* \asymp \left(\frac{1}{n}\right)^{\frac{3}{2d+11}},
$$
and we have
$$
\biggl(
 AMSE(b_1^*)
 +
 R_n(b_0^*, b_1^*)
\biggr)^{\frac{1}{2}} \asymp
\left(\frac{1}{n}\right)^{\frac{6}{2d+11}}.
$$
\end{enumerate}
\label{Optimbandw2}
\end{thm}

\vskip 0.3cm\noindent The results of Theorem \ref{Optimbandw2}
show that the rate $n^{-2/5}$ is reachable if and only when $d\leq
2$. These results are derived  from Theorem \ref{Optimbandw1}.
This latter indicates that if $b_1$ is proportional to $n^{-1/5}$,
the bandwidth $b_0^*$ has the same order as
$$
\max \left\lbrace \left(\frac{1}{n}\right)^{\frac{7}{5(d+4)}} ,
\left(\frac{1}{n}\right)^{\frac{8}{5(2d+4)}} \right\rbrace =
\left(\frac{1}{n}\right)^{\frac{8}{5(2d+4)}}.
$$
 For $d\leq 2$, this order of $b_0^*$ is smaller
than the one of the optimal bandwidth $b_{0*}$ obtained  for
pointwise or mean square estimation of $m(\cdot)$ using a Kernel
estimator. In fact, it has been shown in  Nadaraya (1989, Chapter
4) that the optimal bandwidth $b_{0*}$ for estimating $m(\cdot)$
is obtained by minimizing the order of  the risk function
$$
r_n(b_0) = \esp\left[ \int \mathds{1} \left(x\in\mathcal{X}\right)
\left(\widehat{m}_n(x)-m(x)\right)^2
 \widehat{g}_n^2(x)w(x)
 dx
 \right],
$$
where $\widehat{g}_n(x)$ is a nonparametric Kernel estimator of
$g(x)$, and $w(\cdot)$ is a nonnegative weight function, which is
bounded and squared integrable on $\mathcal{X}$. If $g(\cdot)$ and
$m(\cdot)$ have  continuous second order partial derivatives over
their supports, Nadaraya (1989, Chapter 4) shows that $r_n(b_0)$
has the same order as $b_0^4+\left(1/(nb_0^d)\right)$, leading to
the optimal bandwidth $\widehat{b}_{0}=n^{-1/(d+4)}$  for the
convergence of the estimator $\widehat{m}_n(\cdot)$ of $m(\cdot)$
in the set of the square integrable functions on $\mathcal{X}$.

\vskip 0.2cm\noindent For d=1, the optimal order of $b_0^*$ is
$n^{-(1/5)\times(4/3)}$ which goes to 0 slightly faster than
$n^{-1/5}$, the optimal order of the bandwidth $\widehat{b}_{0}$
for the mean square nonparametric estimation of $m(\cdot)$.

\vskip 0.2cm\noindent For $d=2$, the optimal order of $b_0^*$ is
$n^{-1/5}$. Again this order goes to 0 faster than the order
$n^{-1/6}$ of the optimal bandwidth for the nonparametric
estimation of the regression function with two covariates.

\vskip 0.2cm\noindent However, for $d\geq 3$, we note that the
order of $b_0^*$ goes to $0$ slowly than $\widehat{b}_{0}$. Hence
our results show that optimal $\widehat{m}_n(\cdot)$ for
estimating $f(\cdot)$ should use a very small bandwidth $b_0$.
This suggests that $\widehat{m}_n(\cdot)$ should be less biased
and should have a higher variance than the optimal Kernel
regression estimator of the estimation setup. Such a finding
parallels Wang, Cai, Brown and Levine (2008) who show that a
similar result hold when estimating the conditional variance of a
heteroscedastic regression error term. However Wang et {\it al.}
(2008) do not give the order of the optimal bandwidth to be used
for estimating the regression function in their heteroscedastic
setup. These results show that estimators of $m(\cdot)$ with
smaller bias should be preferred in our framework, compared to the
case where the regression function $m(\cdot)$ is the parameter of
interest.

\subsection{Asymptotic normality}
 We give  now  an asymptotic normality of the
estimator $\widehat{f}_{1n}(\epsilon)$.

\begin{thm}

Assume that
$$
{(\rm\bf{A}_{11}):} \quad nb_0^{d+4}=O(1), \quad nb_0^4b_1=o(1),
\quad
 nb_0^{d}b_1^3\rightarrow\infty,
$$
when $n$ goes to $\infty$. Then under $(A_1)-(A_5)$,
$(A_7)-(A_{10})$, we have
$$
\sqrt{nb_1} \left(
\widehat{f}_{1n}(\epsilon)-\overline{f}_{1n}(\epsilon) \right)
\stackrel{d}{\rightarrow} \mathcal{N} \left(
 0,
\frac{f(\epsilon)} {\prob\left(X\in\mathcal{X}_0\right)} \int
K_1^2(v) dv \right),
$$
where
$$
\overline{f}_{1n}(\epsilon)
 =
f(\epsilon) + \frac{b_1^2}{2}f^{(2)}(\epsilon) \int v^2 K_1(v) dv
+ o\left(b_1^2\right).
$$
\label{normalite}
\end{thm}

\noindent The result of this  theorem shows that the best choice
$b_1^*$ for the bandwidth $b_1$ should achieve the minimum of the
Asymptotic Mean Integrated Square Error
$$
{\rm AMISE}
 =
 \frac{b_1^4}{4}
 \int
 \left(f^{(2)}(\epsilon)\right)^2
 d\epsilon
 \left(
 \int
 v^2K_1(v)dv
 \right)^2
 +
 \frac{1}{nb_1\prob\left(X\in\mathcal{X}_0\right)}
 \int K_1^2(v)dv,
$$
leading to the  optimal bandwidth
$$
b_1^* = \left[ \frac{\frac{1}
 {\prob\left(X\in\mathcal{X}_0\right)}
\displaystyle{\int}K_1^2(v)dv} {\displaystyle{\int}
(f^{(2)}(\epsilon))^2d\epsilon \left(
 \displaystyle{\int}
 v^2K_1(v)dv\right)^2}
 \right]^{1/5}
 n^{-1/5}.
$$
We also note that for $d\leq 2$, $b_1=b_1^*$ and $b_0=b_0^*$,
Theorems \ref{Optimbandw2} and  \ref{Optimbandw1} give
$$
b_1 \asymp \left(\frac{1}{n}\right)^{\frac{1}{5}}, \quad b_0
\asymp \left(\frac{1}{n}\right)^{\frac{8}{5(2d+4)}},
$$
which yields that
$$
nb_0^{d+4} \asymp
\left(\frac{1}{n}\right)^{\frac{12-2d}{5(2d+4)}}, \quad nb_0^4b_1
\asymp \left(\frac{1}{n}\right)^{\frac{16-8d}{5(2d+4)}}, \quad
nb_0^db_1^3 \asymp
\left(\frac{1}{n}\right)^{\frac{4d-8}{5(2d+4)}}.
$$
This shows that for $d=1$, the $(\rm\bf{A}_{11})$ is realizable
with the optimal bandwidths $b_0^*$ and $b_1^*$. But with these
bandwidths, the last constraint of $(\rm\bf{A}_{11})$ is not
satisfied for $d=2$, since $nb_0^db_1^3$ is bounded when
$n\rightarrow\infty$.

\section{ Conclusion}
 The aim of this chapter was  to study
  the nonparametric Kernel estimation of the
probability density function of the regression error using the
estimated residuals. The difference between the feasible estimator
which uses the estimated residuals and the unfeasible one using
the true residuals are studied. An optimal choice of the
first-step bandwidth used to estimate the residuals is also
proposed. Again, an asymptotic normality of the feasible Kernel
estimator and its rate-optimality are established. One of the
contributions of this paper is the analysis of  the impact of the
estimated residuals on the regression errors p.d.f. Kernel
estimator.

\vskip 0.3cm
 In our setup, the strategy was to use an approach
based on a two-steps procedure which, in a first step, replaces
the unobserved residuals terms by some nonparametric estimators
$\widehat{\varepsilon}_i$. In a second step, the \lg
pseudo-observations\rg  $\widehat{\varepsilon}_i$ are used to
estimate the p.d.f $f(\cdot)$, as if they were the true
$\varepsilon_i$'s. If proceeding so can remedy the curse of
dimensionality, a challenging issue was to measure the impact of
the estimated residuals on the final estimator of $f(\cdot)$ in
the first nonparametric step, and to find the order of the optimal
first-step bandwidth $b_0$. For this choice of $b_0$, our results
indicates that the optimal bandwidth to be used for estimating the
regression function $m(\cdot)$ should be smaller than the optimal
bandwidth for the mean square estimation of $m(\cdot)$. That is to
say, the best estimator $\widehat{m}_n(\cdot)$ of the regression
function $m(\cdot)$ needed for estimating $f(\cdot)$  should have
a lower bias and a higher variance than the optimal Kernel
regression of the estimation setup. With this appropriate choice
of $b_0$, it has been seen that for $d\leq 2$, the nonparametric
estimator $\widehat{f}_{1n}(\epsilon)$ of $f$ can reach the
optimal rate $n^{-2/5}$, which
 corresponds to the exact consistency rate reached for the
Kernel density estimator of real-valued variable. Hence our main
conclusion is that for $d\leq 2$, the  estimator
$\widehat{f}_{1n}(\epsilon)$  used for estimating $f(\epsilon)$ is
not affected by the curse of dimensionality, since there is no
negative effect coming from the estimation of the residuals on the
final estimator of $f(\epsilon)$.

\setcounter{subsection}{0} \setcounter{equation}{0}
\renewcommand{\theequation}{\thesection.\arabic{equation}}

\section{ Proofs section}

\subsection*{Intermediate Lemmas for Proposition
\ref{Prop} and Theorem \ref{thm1}}

\begin{lem}
Define, for $x\in\mathcal{X}_0$,
$$
\widehat{g}_n(x)
 =
 \frac{1}{nb_0^d}
 \sum_{i=1}^n
 K_0
 \left(
 \frac{X_i-x}{b_0}
 \right),
 \quad
 \overline{g}_n(x)
 =
 \esp\left[\widehat{g}_n(x)\right].
$$
Then under $(A_1)-(A_2)$, $(A_4)$, $(A_7)$ and $(A_9)$,
 we have, when $b_0$ goes to $0$,
$$
\sup_{x\in\mathcal{X}_0} \left|\overline{g}_{n}(x)-g(x)\right|
 =
O\left(b_0^{2}\right), \quad \sup_{x\in\mathcal{X}_0}
 \left|
\widehat{g}_{n}(x)
 -
 \overline{g}_n(x)
\right|
 = O_{\prob}
 \left(
 b_0^4 +
 \frac{ \ln n}{nb_0^d}
 \right)^{1/2},
$$
 and
$$
\sup_{x\in\mathcal{X}_0}
 \left|
 \frac{1}{\widehat{g}_{n}(x)}
   -
\frac{1}{g(x)} \right| =
 O_{\prob}
 \left(
 b_0^4
 +
 \frac{\ln n}{nb_0^d}
\right)^{1/2}.
$$
\label{Estig}
\end{lem}

\begin{lem}
Under $(A_1)-(A_4)$, $(A_7)$ and $(A_9)$, we have
$$
\sup_{x\in\mathcal{X}_0} \left| \widehat{m}_n(x)-m(x)
 \right|
 =
O_{\prob}
 \left(
 b_0^4+\frac{\ln n}{nb_0^d}
 \right)^{1/2}.
$$
\label{Estim}
\end{lem}

\begin{lem}
Define for $(x,y)\in\Rit^d\times\Rit$,
\begin{eqnarray*}
f_{n}(\epsilon|x)
 =
 \frac{\frac{1}{nh_0^dh_1}
 \sum_{i=1}^n
 K_0\left(\frac{X_i-x}{h_0}\right)
 K_1\left(\frac{Y_i-m(x)-\epsilon}{h_1}\right)
 }{\frac{1}{nh_0^d}
 \sum_{i=1}^n
 K_0\left(\frac{X_i-x}{h_0}\right)}\;,
\end{eqnarray*}
 Then under $(A_1)-(A_3)$, $(A_6)-(A_9)$,
  we have, when $n$ goes to infinity,
$$
\widetilde{f}_n(\epsilon|x)
 -
f_{n}(\epsilon|x)
 =
o_{\prob} \left( \frac{1}{nh_0^dh_1} \right)^{1/2}.
$$
\label{Reste}
\end{lem}

\begin{lem}
 Set, for $(x,y)\in\Rit^d\times\Rit$,
$$
\widetilde{\varphi}_{in}(x, y)
 =
 \frac{1}{h_0^dh_1}
 K_0\left(\frac{X_i-x}{h_0}\right)
 K_1\left(\frac{Y_i-y}{h_1} \right).
$$
Then, under  $(A_6)-(A_8)$, we have, for $x$ in
 $\mathcal{X}_0$ and $y$ in $\Rit$,
  $h_0$ and $h_1$ going to $0$, and for some constant $C>0$,
\begin{eqnarray*}
\esp \left[ \widetilde{\varphi}_{in} \left(x, y\right) \right] -
\varphi\left(x, y\right) &=& \frac{h_0^2}{2}
\frac{\partial^2\varphi (x, y)} {\partial^2 x} \int z
K_0(z)z^{\top} dz
 +
\frac{h_1^2}{2} \frac{\partial^2\varphi(x, y)} {\partial^2 y} \int
v^2 K_1(v) dv
\\
&& +\;
 o\left(h_0^{2}+h_1^2\right),
\\
\Var \left[ \widetilde{\varphi}_{in} \left(x, y\right) \right]
 &=&
\frac{\varphi\left(x, y\right)} {h_0^dh_1}
 \int\int
 K_0^2(z)
 K_1^2(v)
 dv dz
+ o\left(\frac{1}{h_0^dh_1} \right),
\\
\esp \left[ \left| \widetilde{\varphi}_{in}\left(x,y\right) - \esp
\widetilde{\varphi}_{in}\left(x, y\right) \right|^3 \right] &\leq&
 \frac
{C\varphi\left(x,y\right)} {h_0^{2d}h_1^2} \int \int \left| K_0(z)
K_1\left(v\right) \right|^3 dz dv +
o\left(\frac{1}{h_0^{2d}h_1^2}\right).
\end{eqnarray*}
\label{Momvarphi}
\end{lem}

\begin{lem}
Set
$$
f_{in}(\epsilon) = \frac{\mathds{1}\left(X_i \in
\mathcal{X}_0\right)} {b_1\prob\left(X\in\mathcal{X}_0\right)}
K_1\left( \frac{\varepsilon_i - \epsilon}{b_1}\right).
$$
Then under  $(A_4)$, $(A_5)$ and $(A_8)$, we have, for $b_1$ going
to $0$, and for some constant $C>0$,
\begin{eqnarray*}
\esp
 f_{in}(\epsilon)
&=& f(\epsilon) + \frac{b_1^2}{2} f^{(2)}(\epsilon) \int v^2
K_1(v) dv + o\left(b_1^2\right),
\\
\Var \left(f_{in}(\epsilon)\right) &=& \frac {f(\epsilon)}
{b_1\prob\left(X\in\mathcal{X}_0\right)} \int K_1^2(v) dv +
o\left(\frac{1}{b_1}\right),
\\
\esp \left| f_{in}(\epsilon) - \esp
 f_{in}(\epsilon)
\right|^3 &\leq&
 \frac{
C f(\epsilon)}
 {b_1^2\prob^2\left(X\in\mathcal{X}_0\right)}
 \int
 \left|
 K_1(v)
\right|^3 dv + o\left(\frac{1}{b_1^2}\right).
\end{eqnarray*}
\label{Momfin1}
\end{lem}

\begin{lem}
Define
\begin{eqnarray*}
S_n
 &=&
 \sum_{i=1}^n
 \mathds{1}
 \left(X_{i} \in \mathcal{X}_0\right)
 \left(\widehat{m}_{in}-m(X_{i})\right)
 K_1^{(1)}
 \left(\frac{\varepsilon_{i}-\epsilon}{b_1}\right),
 \\
 T_n
 &=&
 \sum_{i=1}^n
 \mathds{1}
 \left(X_{i} \in \mathcal{X}_0\right)
 \left(\widehat{m}_{in}-m(X_{i})\right)^2
 K_1^{(2)}
 \left(\frac{\varepsilon_{i}-\epsilon}{b_1}\right),
 \\
 R_n
 &=&
 \sum_{i=1}^n
 \mathds{1}
 \left( X_{i}
 \in\mathcal{X}_0\right)
 \left(\widehat{m}_{in}-m(X_{i})\right)^3
 \int_{0}^{1}
 (1-t)^2
 K_1^{(3)}
 \left(
 \frac{
 \varepsilon_i-t(\widehat{m}_{in} -m(X_i))-\epsilon}{b_1}
 \right)
 dt.
\end{eqnarray*}
Then under $(A_1)-(A_5)$ and $(A_7)-(A_{10})$, we have, for $b_0$
and $b_1$ small enough,
\begin{eqnarray*}
S_n
 &=&
 O_{\prob}
\left[ b_0^2\left(nb_1^2+(nb_1)^{1/2}\right)
 +
\left( nb_1^4+\frac{b_1}{b_0^d} \right)^{1/2}
 \right],
 \\
 T_{n}
 &=&
 O_{\prob}
\left[ \left( nb_1^3 + \left(nb_1\right)^{1/2} +
\left(n^2b_0^db_1^{3}\right)^{1/2} \right)
\left(b_0^4+\frac{1}{nb_0^d}\right) \right],
\\
R_n &=& O_{\prob} \left[ \left( nb_1^3 +
\left(n^2b_0^db_1\right)^{1/2} \right)
\left(b_0^4+\frac{1}{nb_0^d}\right)^{3/2} \right].
\end{eqnarray*}
\label{STR}
\end{lem}

\begin{lem}
Under  $(A_5)$ and $(A_8)$ we have, for some constant $C>0$, and
for any $\epsilon$ in $\Rit$ and $p\in[0,2]$,
\begin{eqnarray}
\left|
 \int
 K_1^{(1)}
 \left( \frac{e-\epsilon}{b_1} \right)^2
 e^pf(e) de
 \right|
 \leq C b_1,
 &&
 \left|
  \int
  K_1^{(1)}
  \left(
 \frac{e-\epsilon}{b_1}
 \right)
 e^pf(e)de
  \right|
 \leq C b_1^2,
\label{MomderK1}
\\
\left| \int
 K_1^{(2)}
 \left( \frac{e-\epsilon}{b_1} \right)^2
  e^pf(e)
  de
\right|
 \leq
  C b_1,
  &&
 \left|
\int K_1^{(2)}
 \left(
 \frac{e-\epsilon}{b_1}
 \right)
 e^pf(e)de
\right|
 \leq C b_1^3,
 \label{MomderK2}
\\
\left| \int
 K_1^{(3)}
 \left( \frac{e-\epsilon}{b_1} \right)^2
  e^pf(e)
  de
\right|
 \leq
  C b_1,
  &&
 \left|
\int K_1^{(3)}
 \left(
 \frac{e-\epsilon}{b_1}
 \right)
 e^pf(e)de
\right|
 \leq
 C b_1^3.
 \label{MomderK3}
\end{eqnarray}
\label{MomderK}
\end{lem}

\begin{lem}
Set
$$
\beta_{in}
 =
\frac{\mathds{1} \left(X_i\in\mathcal{X}_0\right)}
{nb_0^d\widehat{g}_{in}} \sum_{j=1, j\neq i}^n
\left(m(X_j)-m(X_i)\right) K_0\left(\frac{X_j-X_i}{b_0}\right).
$$
Then, under $(A_1)-(A_5)$ and $(A_7)-(A_{10})$, we have, when
$b_0$ and $b_1$ go to $0$,
\begin{eqnarray*}
\sum_{i=1}^n \beta_{in}
 K_1^{(1)}
 \left(
 \frac{\varepsilon_i-\epsilon}{b_1}
 \right)
 =
 O_{\prob}
 \left(b_0^2\right)
 \left(nb_1^2+(nb_1)^{1/2}\right).
\end{eqnarray*}
\label{Betasum}
\end{lem}

\begin{lem}
Set
$$
\Sigma_{in}
 =
\frac{\mathds{1}\left(X_i\in\mathcal{X}_0\right)}
{nb_0^d\widehat{g}_{in}} \sum_{j=1, j\neq i}^n
 \varepsilon_j
K_0\left(\frac{X_j-X_i}{b_0}\right).
$$
Then, under $(A_1)-(A_5)$ and $(A_7)-(A_{10})$, we have
\begin{eqnarray*}
\sum_{i=1}^n \Sigma_{in}
 K_1^{(1)}
 \left(
 \frac{\varepsilon_i-\epsilon}{b_1}
 \right)
 =
 O_{\prob}
 \left(
 nb_1^4
 +
 \frac{b_1}{b_0^d}
\right)^{1/2}.
\end{eqnarray*}
\label{Sigsum}
\end{lem}

\begin{lem}
Let $\esp_n[\cdot]$ be the conditional mean given
$X_1,\ldots,X_n$. Then under $(A_1)-(A_5)$ and $(A_7)-(A_{9})$,
 we have, for $b_0$ going to $0$,
\begin{eqnarray*}
\sup_{1\leq i\leq n}
 \esp_{n}
  \biggl[
  \mathds{1}
  \left(X_i \in \mathcal{X}_0\right)
 (\widehat{m}_{in} - m(X_i))^4
 \biggr]
 &=&
 O_{\prob}
 \left(b_0^4+\frac{1}{nb_0^d}\right)^2,
 \\
 \sup_{1\leq i\leq n}
 \esp_{n}
 \biggl[
 \mathds{1}
 \left(X_i \in \mathcal{X}_0\right)
 (\widehat{m}_{in} - m(X_i))^6
 \biggr]
 &=&
 O_{\prob}
 \left(b_0^4+\frac{1}{nb_0^d}\right)^3.
\end{eqnarray*}
\label{BoundEspmchap}
\end{lem}

\begin{lem}
Assume that  $(A_4)$ and $(A_7)$ hold. Then, for any $1 \leq i
\neq j \leq n$, and for any $\epsilon$ in $\Rit$,
$$
\left(\widehat{m}_{in} - m(X_i), \varepsilon_i\right)
 \mbox{\it and }
 \left( \widehat{m}_{jn} - m(X_j), \varepsilon_j\right)
$$
are independent given $X_1, \ldots, X_n$,
 provided that $\| X_i - X_j \| \geq C b_0$,
for some constant $C>0$. \label{Indep}
\end{lem}

\begin{lem}
Let $\Var_n(\cdot)$ and $Cov_n(\cdot)$ be respectively the
conditional variance and the conditional covariance given
$X_1,\ldots,X_n$, and set
\begin{eqnarray*}
\zeta_{in}
 =
\mathds{1} \left(X_i \in \mathcal{X}_0\right) (\widehat{m}_{in} -
m(X_i))^2 K_1^{(2)} \left( \frac{\varepsilon_i-\epsilon}{b_1}
\right).
\end{eqnarray*}
Then under $(A_1)-(A_5)$ and $(A_7)-(A_{9})$, we have, for $n$
going to infinity,
\begin{eqnarray*}
\sum_{i=1}^n \Var_n\left( \zeta_{in}\right) &=&
O_{\prob}\left(nb_1\right)
 \left(b_0^4+\frac{1}{nb_0^d}\right)^2,
 \\
 \sum_{i=1}^n
 \sum_{j=1\atop j\neq i}^n
  \Cov_n
 \left(
 \zeta_{in},\zeta_{jn}
 \right)
&=& O_{\prob}\left(n^2b_0^db_1^{7/2}\right)
\left(b_0^4+\frac{1}{nb_0^d}\right)^2.
\end{eqnarray*}
\label{sumzeta}
\end{lem}

\noindent All these lemmas are proved in Appendix A.

\subsection*{Proof of Proposition \ref{Prop}}
Define $f_{n}(\epsilon|x)$ as in Lemma \ref{Reste}, and note that
by this lemma, we have
 \begin{eqnarray}
 \widetilde{f}_n(\epsilon|x)
 =
 f_{n}(\epsilon|x)
 +
 o_{\prob}
\left( \frac{1}{nh_0^dh_1} \right)^{1/2}.
 \label{fncond1}
\end{eqnarray}
The asymptotic distribution of the first term in (\ref{fncond1})
is derived by applying the Lyapounov Central Limit Theorem for
triangular arrays (see e.g Billingsley 1968, Theorem 7.3). Define
for $x\in\mathcal{X}_0$ and $y\in\Rit$,
\begin{eqnarray*}
 \widetilde{\varphi}_n(x,y)
 =
 \frac{1}{nh_0^dh_1}
\sum_{i=1}^n
 K_0\left(\frac{X_i-x}{h_0}\right)
 K_1\left(\frac{Y_i-y}{h_1}\right),
 \quad
 \widetilde{g}_n(x)
 =
 \frac{1}{nh_0^d}
 \sum_{i=1}^n
 K_0\left(\frac{X_i-x}{h_0}\right),
\end{eqnarray*}
and observe that
\begin{eqnarray}
f_{n}(\epsilon|x)
 =
\frac{\widetilde{\varphi}_n\left(x,m(x)+\epsilon\right)}
{\widetilde{g}_n(x)}\;.
 \label{fncond2}
\end{eqnarray}
Let now $\widetilde{\varphi}_{in}(x,y)$ be as in Lemma
\ref{Momvarphi}, and note that
\begin{eqnarray}
\widetilde{\varphi}_n(x,y)
 =
 \frac{1}{n}
 \sum_{i=1}^n
 \biggl(
 \widetilde{\varphi}_{in}(x,y)
 -
 \esp
 \left[
 \widetilde{\varphi}_{in}(x,y)
 \right]
 \biggr)
 +
 \esp
 \left[
 \widetilde{\varphi}_{1n}(x,y)
 \right].
 \label{fnoncond}
\end{eqnarray}
The second and third inequalities in Lemma \ref{Momvarphi} give,
since  $h_0^d h_1$ goes to $0$,
\begin{eqnarray*}
\frac{ \sum_{i=1}^n \esp \left|
 \widetilde{\varphi}_{in}(x,y)
 -
 \esp
 \widetilde{\varphi}_{in}(x,y)
\right|^3 } { \left( \sum_{i=1}^n
 \Var
 \left[
 \widetilde{\varphi}_{in}(x,y)
 \right]
\right)^3 } \leq \frac { \frac { Cn\varphi(x,y) } {h_0^{2d}h_1^2}
\displaystyle{\int} \displaystyle{\int} \left| K_0(z) K_1(v)
\right|^3 dz dv + o \left( \frac{n}{h_0^{2d}h_1^2} \right) } {
\left( \frac { n\varphi(x,y) }{h_0^{d}h_1}
 \displaystyle{\int}
\displaystyle{\int}
 K_0^2(z) K_1^2(v) dv dz + o \left(
\frac{n}{h_0^{d}h_1} \right)
 \right)^3
 }
 =O(h_0^d h_1) = o(1).
\end{eqnarray*}
 Hence  the Lyapounov Central Limit
Theorem gives, since $nh_0^dh_1$ diverges under $(\rm\bf{A}_0)$,
$$
\frac { \sum_{i=1}^n
 \left\lbrace
 \widetilde{\varphi}_{in}(x,y)
 -
\esp \left[ \widetilde{\varphi}_{in}(x,y) \right] \right\rbrace}{
\left( \sum_{i=1}^n \Var \left[ \widetilde{\varphi}_{in}(x,y)
\right] \right)^{1/2}} \stackrel{d}{\rightarrow}\mathcal{N}(0,1),
$$
so that
\begin{eqnarray}
\frac{\sqrt{nh_0^dh_1}}{n}
 \sum_{i=1}^n
 \biggl(
 \widetilde{\varphi}_{in}(x,y)
 -
 \esp
 \left[
 \widetilde{\varphi}_{in}(x,y)
 \right]
 \biggr)
 \stackrel{d}{\rightarrow}
 \mathcal{N}
 \left(
  0,
  \varphi(x,y)
  \int\int
  K_0^2(z)
  K_1^2(v)
  dz dv
 \right).
\label{normal}
\end{eqnarray}
Further, a similar proof as the one of Lemma \ref{Estig} gives
\begin{eqnarray}
\frac{1}{\widetilde{g}_n(x)}
 =
 \frac{1}{g(x)}
 +
 O_{\prob}
 \left(
 h_0^4+\frac{\ln n}{nh_0^d}
 \right)^{1/2}.
 \label{gtilde}
\end{eqnarray}
Hence by this equality, it follows that, taking $y=m(x)+\epsilon$
in (\ref{normal}), and by (\ref{fncond1})-(\ref{fnoncond}),
\begin{eqnarray*}
\sqrt{nh_0^dh_1}
 \left(
\widetilde{f}_n(\epsilon|x)
 -
\overline{f}_n(\epsilon|x) \right) \stackrel{d}{\rightarrow}
\mathcal{N} \left( 0, \frac{f(\epsilon|x)}{g(x)} \int\int
 K_0^2(z)
 K_1^2(v)
 dz dv
\right),
\end{eqnarray*}
where
$$
\overline{f}_n(\epsilon|x) = \frac{ \esp\left[
\widetilde{\varphi}_{1n} \left(x,m(x)+\epsilon\right) \right]}
{\widetilde{g}_n(x)}\;.
$$
This yields  the result of Proposition \ref{Prop}, since the first
equality of Lemma \ref{Momvarphi} and (\ref{gtilde}) yield, for
$h_0$ and $h_1$ small enough,
\begin{eqnarray*}
\overline{f}_n(\epsilon|x)
 &=&
 f(\epsilon| x)
+ \frac{h_0^2}{2g(x)} \frac{\partial^2\varphi\left(x,
m(x)+\epsilon\right)} {\partial^2 x} \int z K_0(z) z^{\top} dz
 \\
 &&
 +\;
 \frac{h_1^2}{2g(x)}
 \frac{\partial^2\varphi (x,m(x)+\epsilon)}
 {\partial^2 y}
 \int v^2 K_1 (v) dv
 +
 o\left(h_0^2+h_1^2\right).
 \eop
\end{eqnarray*}

\subsection*{Proof of Theorem \ref{thm1}}
The proof of the theorem is based upon the following equalities:
\begin{eqnarray}
\nonumber
 \widehat{f}_{1n}(\epsilon)-\widetilde{f}_{1n}(\epsilon)
&=& O_{\prob}
 \left[
 b_0^2
 +
 \left(
 \frac{1}{n}
 +
 \frac{1}{n^2b_0^db_1^3}
 \right)^{1/2}
 \right]
 +
 O_{\prob}
 \left[
 \frac{1}{(nb_1^5)^{1/2}}
 +
 \left(\frac{b_0^d}{b_1^3}\right)^{1/2}
 \right]
 \left(
 b_0^4
 +
 \frac{1}{nb_0^d}
 \right)
 \\
 &&
 +\;
 O_{\prob}
 \left[
 \frac{1}{b_1}
 +
 \left(\frac{b_0^d}{b_1^7}\right)^{1/2}
 \right]
 \left(b_0^4+\frac{1}{nb_0^d}\right)^{3/2},
 \label{fnchapfn1}
 \end{eqnarray}
 and
\begin{eqnarray}
\widetilde{f}_{1n}(\epsilon)-f(\epsilon)
 =
 O_{\prob}
 \left(b_1^4+\frac{1}{nb_1}\right)^{1/2}.
 \label{fnchapfn2}
 \end{eqnarray}
Indeed, since $\widehat{f}_{1n}(\epsilon)-f(\epsilon)
=\left(\widetilde{f}_{1n}(\epsilon)-f(\epsilon)\right)
+\widehat{f}_{1n}(\epsilon)-\widetilde{f}_{1n}(\epsilon)$, it then
follows by (\ref{fnchapfn2}) and (\ref{fnchapfn1}) that
\begin{eqnarray*}
\widehat{f}_{1n}(\epsilon)-f(\epsilon) &=& O_{\prob} \left[
 b_1^4
 +
 \frac{1}{nb_1}
 +
 b_0^4
 +
 \frac{1}{n}
 +
 \frac{1}{n^2b_0^db_1^3}
 +
 \left(
 \frac{1}{(nb_1^5)^{1/2}}
 +
 \left(\frac{b_0^d}{b_1^3}\right)^{1/2}
 \right)^2
 \left(
 b_0^4
 +
 \frac{1}{nb_0^d}
 \right)^2
 \right]^{1/2}
 \\
 &&
 +\;
 O_{\prob}
 \left[
 \left(
 \frac{1}{b_1}
 +
 \left(\frac{b_0^d}{b_1^7}\right)^{1/2}
 \right)^2
 \left(b_0^4+\frac{1}{nb_0^d}\right)^3
 \right]^{1/2}.
 \end{eqnarray*}
This yields the result of the Theorem,  since  under $(A_9)$ and
$(A_{10})$, we have
\begin{eqnarray*}
 \frac{1}{n}
  =
 O\left(
 \frac{1}{nb_1}
 \right),
 \quad
\frac{1}{n^2b_0^db_1^3} = O\left(\frac{b_0^d}{b_1^3}\right)
\left(b_0^4+\frac{1}{nb_0^d}\right)^2.
 \end{eqnarray*}
 Hence, it remains
to prove (\ref{fnchapfn1}) and (\ref{fnchapfn2}). For this,
 define $S_n$, $R_n$ and $T_n$ as in Lemma \ref{STR}.
Since $\widehat{\varepsilon}_i-\varepsilon_i =
-\left(\widehat{m}_{in}-m(X_{i})\right)$ and that $K_1$ is three
times continuously differentiable under $(A_8)$,  the third-order
Taylor expansion
 with integral remainder gives
\begin{eqnarray*}
 \widehat{f}_{1n} (\epsilon)- \widetilde{f}_{1n}(\epsilon)
 &=&
\frac{1}{b_1\sum_{i=1}^n \mathds{1}(X_i \in\mathcal{X}_0)}
\sum_{i=1}^n \mathds{1} \left(X_i\in \mathcal{X}_0\right) \left[
 K_1\left(
 \frac{\widehat{\varepsilon}_i-\epsilon}{b_1}
 \right)
 -
 K_1\left(
 \frac{\varepsilon_i-\epsilon}{b_1}
 \right)
 \right]
 \\
& = &
 -\frac{1}{b_1\sum_{i=1}^n\mathds{1}(X_i \in\mathcal{X}_0)}
\left( \frac{S_n}{b_1} - \frac{T_n}{2b_1^2} + \frac{R_n}{2b_1^3}
\right).
\end{eqnarray*}
Therefore, since
$$
\sum_{i=1}^n \mathds{1} \left(X_i \in \mathcal{X}_0\right)
 =
n\left( \prob \left( X \in \mathcal{X}_0\right)
 +
 o_{\prob}(1)
 \right),
$$
 by the Law of large numbers, Lemma \ref{STR} then gives
\begin{eqnarray*}
\lefteqn{ \widehat{f}_{1n}(\epsilon)-\widetilde{f}_{1n}(\epsilon)
 =
O_{\prob}\left(\frac{1}{nb_1^2}\right)S_n + O_{\prob}
\left(\frac{1}{nb_1^3}\right) T_n
 +
O_{\prob} \left(\frac{1}{nb_1^4}\right) R_n }
\\
&=&
 O_{\prob}
 \left[
 b_0^2
 \left(1+\frac{1}{(nb_1^3)^{1/2}}\right)
 +
 \left(
 \frac{1}{n}
 +
 \frac{1}{n^2b_0^db_1^3}
 \right)^{1/2}
 \right]
 \\
 &&
 +\;
 O_{\prob}
 \left[
 1
 +
 \frac{1}{(nb_1^5)^{1/2}}
 +
 \left(\frac{b_0^d}{b_1^3}\right)^{1/2}
 \right]
 \left(
 b_0^4
 +
 \frac{1}{nb_0^d}
 \right)
 +
 O_{\prob}
 \left[
 \frac{1}{b_1}
 +
 \left(\frac{b_0^d}{b_1^7}\right)^{1/2}
 \right]
 \left(b_0^4+\frac{1}{nb_0^d}\right)^{3/2}.
 \end{eqnarray*}
This yields (\ref{fnchapfn1}), since  under $(A_9)$ and
$(A_{10})$,  we have $b_0\rightarrow 0$,
$nb_0^{d+2}\rightarrow\infty$ and  $nb_1^3\rightarrow\infty$, so
that
\begin{eqnarray*}
b_0^2
 \left(1+\frac{1}{(nb_1^3)^{1/2}}\right)
&\asymp&
 O\left(b_0^2\right),
\quad
 \left(
 b_0^4
 +
 \frac{1}{nb_0^d}
 \right)
 =
 O\left(b_0^2\right),
\\
\left[
 1
 +
 \frac{1}{(nb_1^5)^{1/2}}
 +
 \left(\frac{b_0^d}{b_1^3}\right)^{1/2}
 \right]
 \left(
 b_0^4
 +
 \frac{1}{nb_0^d}
 \right)
 &=&
 O\left(b_0^2\right)
 +
 \left[
 \frac{1}{(nb_1^5)^{1/2}}
 +
 \left(\frac{b_0^d}{b_1^3}\right)^{1/2}
 \right]
 \left(
 b_0^4
 +
 \frac{1}{nb_0^d}
 \right).
\end{eqnarray*}

\vskip 0.1cm
 For (\ref{fnchapfn2}), note that
\begin{eqnarray}
\esp_n \left[
\left(\widetilde{f}_{1n}(\epsilon)-f(\epsilon)\right)^2 \right] =
\Var_n\left(\widetilde{f}_{1n}(\epsilon)\right) + \biggl( \esp_n
\left[ \widetilde{f}_{1n}(\epsilon)\right] - f(\epsilon)
\biggr)^2, \label{Quadbias}
\end{eqnarray}
with, using $(A_4)$,
\begin{eqnarray*}
\Var_n\left(\widetilde{f}_{1n}(\epsilon)\right)
 =
\frac{1}{\left(b_1\sum_{i=1}^n
\mathds{1}\left(X_i\in\mathcal{X}_0\right)\right)^2} \sum_{i=1}^n
\mathds{1} \left(X_i\in\mathcal{X}_0\right) \Var \left[
 K_1\left(\frac{\varepsilon-\epsilon}{b_1}\right)
 \right].
\end{eqnarray*}
Therefore, since the Cauchy-Schwarz inequality gives
\begin{eqnarray*}
\Var \left[
 K_1\left(\frac{\varepsilon-\epsilon}{b_1}\right)
 \right]
\leq
 \esp
 \left[
 K_1^2\left(\frac{\varepsilon-\epsilon}{b_1}\right)
 \right]
 \leq
 b_1\int
 K_1^2(v)f(\epsilon+b_1v)
 dv,
\end{eqnarray*}
this bound and the equality above yield, under $(A_5)$ and
$(A_8)$,
\begin{eqnarray}
 \Var_n\left(\widetilde{f}_{1n}(\epsilon)\right)
 \leq
 \frac{C}{b_1\sum_{i=1}^n
 \mathds{1}\left(X_i\in\mathcal{X}_0\right)}
 =
 O_{\prob}\left(\frac{1}{nb_1}\right).
 \label{Quadbias1}
\end{eqnarray}
For the second term in (\ref{Quadbias}),  we have
\begin{eqnarray}
\esp_n\left[\widetilde{f}_{1n}(\epsilon)\right] =
\frac{1}{b_1\sum_{i=1}^n
 \mathds{1}\left(X_i\in\mathcal{X}_0\right)}
\sum_{i=1}^n \mathds{1} \left(X_i\in\mathcal{X}_0\right) \esp
\left[
 K_1\left(\frac{\varepsilon-\epsilon}{b_1}\right)
 \right].
 \label{Quadbias2}
\end{eqnarray}
By $(A_8)$,  $K_1$  is symmetric, has  a compact support, with
 $\int\! v K_1(v)=0$  and $\int \! K_1(v)dv=1$.
 Therefore, since under $(A_5)$ $f$ has bounded continuous second
order derivatives, this yields for some
$\theta=\theta(\epsilon,b_1v)$,
\begin{eqnarray*}
\lefteqn{
 \esp
 \left[
 K_1\left(\frac{\varepsilon-\epsilon}{b_1}\right)
 \right]
 =
  b_1
\int
 K_1(v)
 f(\epsilon+b_1v)
 dv
 }
\\
&=&
 b_1
 \int
 K_1(v)
\left[
 f(\epsilon)
+
 b_1 v f^{(1)}(\epsilon)
 +
 \frac{b_1^2v^2}{2}
 f^{(2)}(\epsilon+\theta b_1v)
 \right]
 dv
 \\
  &=&
 b_1f(\epsilon)
 +
 \frac{b_1^3}{2}
 \int v^2 K_1(v)
 f^{(2)}(\epsilon+\theta b_1v)
 dv.
\end{eqnarray*}
Hence  this equality and (\ref{Quadbias2}) give
$$
\esp_n \left[\widetilde{f}_{1n}(\epsilon)\right]
 =
 f(\epsilon)
 +
 \frac{b_1^2}{2}
 \int v^2 K_1(v)
 f^{(2)}(\epsilon+\theta b_1v)
 dv,
$$
so that
$$
\biggl( \esp_n \left[\widetilde{f}_{1n}(\epsilon)\right] -
f(\epsilon) \biggr)^2
 =
 O_{\prob}\left(b_1^4\right).
$$
Combining this result with  (\ref{Quadbias1}) and
(\ref{Quadbias}), we obtain, by the Tchebychev inequality,
$$
\widetilde{f}_{1n}(\epsilon)-f(\epsilon) = O_{\prob}
\left(b_1^4+\frac{1}{nb_1}\right)^{1/2}.
$$
This proves (\ref{fnchapfn2}), and then achieves the proof of the
theorem.\eop

\subsection*{Proof of Theorem \ref{Optimbandw1}}
Recall that
\begin{eqnarray*}
 R_n(b_0, b_1)
=
 b_0^4
 +
 \left[
 \frac{1}{(nb_1^5)^{1/2}}
 +
 \left(\frac{b_0^d}{b_1^3}\right)^{1/2}
 \right]^2
 \left(
 b_0^4
 +
 \frac{1}{nb_0^d}
 \right)^2
 +
 \left[
 \frac{1}{b_1}
 +
 \left(\frac{b_0^d}{b_1^7}\right)^{1/2}
 \right]^2
 \left(b_0^4+\frac{1}{nb_0^d}\right)^3,
\end{eqnarray*}
and note that
$$
\left(\frac{1}{n^2b_1^3}\right)^{\frac{1}{d+4}} = \max
\left\lbrace
 \left(\frac{1}{n^2b_1^3}\right)^{\frac{1}{d+4}}
, \left(\frac{1}{n^3b_1^7}\right)^{\frac{1}{2d+4}} \right\rbrace
$$
if and only if $n^{4-d}b_1^{d+16}\rightarrow\infty$. To find the
order of $b_0^*$, we shall deal with the cases
$nb_0^{d+4}\rightarrow\infty$ and
 $nb_0^{d+4}=O(1)$.
\vskip 0.1cm\noindent First assume that
$nb_0^{d+4}\rightarrow\infty$. More precisely, we  suppose that
$b_0$ is in $\left[(u_n/n)^{1/(d+4)},+\infty\right)$, where $u_n
\rightarrow\infty$. Since $1/(nb_0^d) = O(b_0^4)$ for all these
$b_0$, we have
$$
\left(b_0^4+\frac{1}{nb_0^d}\right)^2 \asymp \left(b_0^4\right)^2,
\quad \left(b_0^4+\frac{1}{nb_0^d}\right)^3 \asymp
 \left(b_0^4\right)^3.
$$
Hence the order of $b_0^*$ is computed  by minimizing the function
\begin{eqnarray*}
 b_0\rightarrow
  b_0^4
 +
 \left[
 \frac{1}{(nb_1^5)^{1/2}}
 +
 \left(\frac{b_0^d}{b_1^3}\right)^{1/2}
 \right]^2
 \left( b_0^4\right)^2
 +
 \left[
 \frac{1}{b_1}
 +
 \left(\frac{b_0^d}{b_1^7}\right)^{1/2}
 \right]^2
 \left(b_0^4\right)^3.
\end{eqnarray*}
Since this function is increasing with $b_0$, the minimum of $R_n
(\cdot,b_1)$ is achieved for $b_{0*}=(u_n/n)^{1/(d+4)}$. We shall
prove later on that this choice of $b_{0*}$ is irrelevant compared
to the one arising when $nb_0^{d+4}=O(1)$.

\vskip 0.1cm Consider now the case $nb_0^{d+4}=O(1)$ i.e
$b_0^4=O\left(1/(nb_0^d)\right)$. This  gives
\begin{eqnarray*}
 \left[
 \frac{1}{(nb_1^5)^{1/2}}
 +
 \left(\frac{b_0^d}{b_1^3}\right)^{1/2}
 \right]^2
 \left(
 b_0^4
 +
 \frac{1}{nb_0^d}
 \right)^2
 &\asymp&
  \left(
 \frac{1}{nb_1^5}
 +
 \frac{b_0^d}{b_1^3}
 \right)
 \left(
 \frac{1}{n^2b_0^{2d}}
 \right),
 \\
 \left[
 \frac{1}{b_1}
 +
 \left(\frac{b_0^d}{b_1^7}\right)^{1/2}
 \right]^2
 \left(b_0^4+\frac{1}{nb_0^d}\right)^3
 &\asymp&
 \left(
 \frac{1}{b_1^2}
 +
 \frac{b_0^d}{b_1^7}
 \right)
 \left(\frac{1}{n^3b_0^{3d}}\right).
 \end{eqnarray*}
Moreover if $nb_0^db_1^4\rightarrow\infty$, we have, since
$nb_0^{2d}\rightarrow\infty$ under $(A_9)$,
$$
 \left(
 \frac{1}{nb_1^5}
 +
 \frac{b_0^d}{b_1^3}
 \right)
 \left(
 \frac{1}{n^2b_0^{2d}}
 \right)
 \asymp
 \frac{b_0^d}{b_1^3}
 \left(
 \frac{1}{n^2b_0^{2d}}
 \right),
 \quad
 \left(
 \frac{1}{b_1^2}
 +
 \frac{b_0^d}{b_1^7}
 \right)
 \left(\frac{1}{n^3b_0^{3d}}\right)
 =
 O\left(
 \frac{b_0^d}{b_1^3}
 \right)
 \left(
 \frac{1}{n^2b_0^{2d}}
 \right).
 $$
Hence the order of $b_0^*$ is obtained by finding the minimum of
the function $b_0^4+\left(1/n^2b_0^db_1^3\right)$. The
 minimization of this function  gives  a solution $b_0$ such that
$$
b_0 \asymp \left(\frac{1}{n^2b_1^3}\right)^{\frac{1}{d+4}}, \quad
R_n(b_0,b_1) \asymp
\left(\frac{1}{n^2b_1^3}\right)^{\frac{4}{d+4}}.
$$
This value  satisfies the constraints $nb_0^{d+4}=O(1)$ and
$nb_0^db_1^4\rightarrow\infty$ when
$n^{4-d}b_1^{d+16}\rightarrow\infty$.

\vskip 0.1cm\noindent If now $nb_0^{d+4}=O(1)$ but
$nb_0^db_1^4=O(1)$, we have, since $nb_0^{2d}\rightarrow\infty$,
\begin{eqnarray*}
\frac{1}{nb_1^5}
 \left(
\frac{1}{n^2b_0^{2d}} \right)
  =
 O\left(
 \frac{b_0^d}{b_1^7}
 \right)
 \left(\frac{1}{n^3b_0^{3d}}\right),
 \quad
 \frac{1}{b_1^2}
 \left(\frac{1}{n^3b_0^{3d}}\right)
 =
O\left(
 \frac{b_0^d}{b_1^3}
 \right)
 \left(
\frac{1}{n^2b_0^{2d}} \right)
 =
 O\left(
 \frac{b_0^d}{b_1^7}
 \right)
 \left(\frac{1}{n^3b_0^{3d}}\right).
\end{eqnarray*}
 In this case, $b_0^*$ is obtained by minimizing the function
$b_0^4+\left(1/n^3b_0^{2d}b_1^7\right)$, for which the solution
 $b_0$ verifies
$$
b_0 \asymp \left(\frac{1}{n^3b_1^7}\right)^{\frac{1}{2d+4}}, \quad
R_n(b_0,b_1) \asymp
\left(\frac{1}{n^3b_1^7}\right)^{\frac{4}{2d+4}}.
$$
This solution fulfills  the constraint $nb_0^db_1^4=O(1)$ when
$n^{4-d}b_1^{d+16}=O(1)$. Hence we can conclude that for
$b_0^4=O\left(1/(nb_0^d)\right)$, the bandwidth $b_0^*$ satisfies
$$
b_0^* \asymp \max \left\lbrace
\left(\frac{1}{n^2b_1^3}\right)^{\frac{1}{d+4}} ,
\left(\frac{1}{n^3b_1^7}\right)^{\frac{1}{2d+4}} \right\rbrace,
$$
which leads to
$$
R_n\left(b_0^*, b_1\right) \asymp \max \left\lbrace
\left(\frac{1}{n^2b_1^3}\right)^{\frac{4}{d+4}} ,
\left(\frac{1}{n^3b_1^7}\right)^{\frac{4}{2d+4}}
 \right\rbrace.
$$
We need now to compare the solution $b_0^*$ to the candidate
 $b_{0*}=(u_n/n)^{1/(d+4)}$ obtained when $nb_0^{d+4}\rightarrow\infty$.
  For this, we must  do a comparison between the
 orders of $R_n(b_0^*,b_1)$  and $R_n(b_{0*},b_1)$. Since
 $R_n(b_0,b_1)\geq b_0^4$, we have
 $R_n(b_{0*},b_1)\geq(u_n/n)^{4/(d+4)}$, so that, for $n$ large
 enough,
\begin{eqnarray*}
\frac{R_n(b_0^*, b_1)} {R_n(b_{0*},b_1)} &\leq&
 C\left[
\left(\frac{1}{n^2b_1^3}\right)^{\frac{1}{d+4}}
 +
\left(\frac{1}{n^3b_1^7}\right)^{\frac{4}{2d+4}} \right]
\left(\frac{n}{u_n}\right)^{\frac{4}{d+4}}
\\
&=& o(1) + O\left(\frac{1}{u_n}\right)^{\frac{4}{d+4}} \left(
\frac{1}{nb_1^{\frac{7(d+4)}{d+8}}}
\right)^{\frac{4(d+8)}{(2d+4)(d+4)}} = o(1),
\end{eqnarray*}
using $u_n\rightarrow\infty$ and that
$n^{(d+8)}b_1^{7(d+4)}\rightarrow\infty$ by $(A_{10})$. This shows
that $R_n(b_0^*,b_1)\leq R_n(b_{0*},b_1)$ for $n$ large enough.
 Hence the Theorem is proved, since $b_0^*$ is the best candidate for
 the minimization of $R_n(\cdot, b_1)$. \eop

\subsection*{Proof of Theorem \ref{Optimbandw2}}
Recall that  Theorem \ref{Optimbandw1} gives
\begin{eqnarray*}
AMSE(b_1)+R_n(b_0^*, b_1) \asymp r_1(b_1)+r_2(b_1)+r_3(b_1) =
F(b_1),
\end{eqnarray*}
where
\begin{eqnarray*}
r_1(h) &=& h^4 + \frac{1}{nh}, \quad \arg\min r_1(h) \asymp
n^{-1/5}=h_1^*, \quad \min r_1(h) \asymp (h_1^*)^4=n^{-4/5},
\\
r_2(h) &=& h^4 + \frac{1}{n^{\frac{8}{d+4}} h^{\frac{12}{d+4}}},
\quad \arg\min r_2(h) \asymp n^{-\frac{2}{d+7}} = h_2^*, \quad
\min r_3(h) \asymp(h_2^*)^4 = n^{-\frac{8}{d+7}},
\\
r_3(h) &=& h^4 + \frac{1}{n^{\frac{12}{2d+4}}
h^{\frac{28}{2d+4}}}, \quad \arg\min r_3(h) \asymp
n^{-\frac{3}{2d+11}} = h_3^*, \quad \min r_3(h) \asymp(h_3^*)^4 =
n^{-\frac{12}{2d+11}}.
\end{eqnarray*}
Each $r_j(h)$ decreases on $\left[0,\arg\min r_j(h)\right]$ and
increases on $\left(\arg\min r_j(h),\infty\right)$ and that
$r_j(h)\asymp h^4$ on $\left(\arg\min r_j(h),\infty\right)$.
Moreover $\min r_2(h)=o\left(r_3(h)\right)$ and
$h_2^*=o\left(h_3^*\right)$ for all possible dimension $d$, so
that $\min\{r_2(h)+r_3(h)\}\asymp(h_3^*)^4=n^{-\frac{12}{2d+11}}$
and $\arg\min\{r_2(h)+r_3(h)\}\asymp h_3^*=n^{-\frac{3}{2d+11}}$.

\vskip 0.1cm Observe now that $\min\{r_2(h)+r_3(h)\}=O\left(\min
r_1(h)\right)$ is equivalent to
$n^{-\frac{12}{2d+11}}=O\left(n^{-4/5}\right)$ which holds if and
only if $d\leq 2$. Hence assume that $d\leq 2$. Since
$n^{-\frac{12}{2d+11}}=O\left(n^{-4/5}\right)$ also gives
$\arg\min\{r_2(h)+r_3(h)\}\asymp h_3^*=O\left(h_1^*\right)$, we
have
$$
\min F(b_1) \asymp n^{-4/5} \;\;{\rm and}\; \arg\min F(b_1) \asymp
n^{-1/5}.
$$
The case $d>2$ is symmetric with
$$
\min F(b_1) \asymp n^{-\frac{12}{2d+11}} \;\;{\rm and}\; \arg\min
F(b_1) \asymp n^{-\frac{3}{2d+11}}.
$$
This ends the proof of the Theorem. \eop

\subsection*{Proof of Theorem \ref{normalite}}
Observe that the Tchebychev inequality gives
$$
\sum_{i=1}^n \mathds{1} \left(X_i\in\mathcal{X}_0\right) =
n\prob\left(X\in\mathcal{X}_0\right) \left[ 1 + O_{\prob}
\left(\frac{1}{\sqrt{n}}\right) \right],
$$
so that
$$
\widetilde{f}_{1n}(\epsilon) = \left[ 1 + O_{\prob}
\left(\frac{1}{\sqrt{n}}\right) \right] f_n(\epsilon),
$$
where
$$
f_n(\epsilon)
 =
\frac{1}{nb_1\prob\left(X\in\mathcal{X}_0\right)} \sum_{i=1}^n
\mathds{1} \left(X_i\in\mathcal{X}_0\right)
 K_1\left(\frac{\varepsilon_i-\epsilon}{b_1}\right).
$$
Therefore
\begin{eqnarray}
\widehat{f}_{1n}(\epsilon)-\esp f_n(\epsilon) = \left(
f_n(\epsilon)-\esp f_n(\epsilon) \right) + \left(
\widehat{f}_{1n}(\epsilon)-\widetilde{f}_{1n}(\epsilon) \right) +
O_{\prob} \left(\frac{1}{\sqrt{n}}\right)
 f_n(\epsilon).
\label{develop}
\end{eqnarray}
Let now $f_{in}(\epsilon)$ be as in Lemma \ref{Momfin1}, and note
that $f_n(\epsilon) =(1/n)\sum_{i=1}^nf_{in}(\epsilon)$. The
second and the third claims in Lemma \ref{Momfin1} yield, since
$b_1$ goes to $0$ under $(A_{10})$,
\begin{eqnarray*}
\frac { \sum_{i=1}^n \esp \left| f_{in}(\epsilon) - \esp
 f_{in}(\epsilon)
\right|^3 } { \left( \sum_{i=1}^n
 \Var f_{in}(\epsilon)
\right)^3 } \leq \frac{ \frac{C nf(\epsilon)}
{\prob\left(X\in\mathcal{X}_0\right)^2b_1^2}
 \displaystyle{\int}
\left| K_1(v) \right|^3 dv + o\left( \frac{n}{b_1^2} \right) } {
\left( \frac{nf(\epsilon)}{\prob\left(
X\in\mathcal{X}_0\right)b_1} \displaystyle{\int}K_1^2(v) dv +
 o\left(\frac{n}{b_1}\right)
\right)^3
 }
 =O(b_1)= o(1).
\end{eqnarray*}
 Hence the Lyapounov Central Limit Theorem gives, since $nb_1$
 diverges under $(A_{10})$,
$$
\frac{f_n(\epsilon)-\esp f_{n}(\epsilon)} {\sqrt{\Var
f_{n}(\epsilon)}} =
 \frac{f_n(\epsilon)-\esp f_{n}(\epsilon)}
 {\sqrt{\frac{\Var f_{in}(\epsilon)}{n}}}
 \stackrel{d}{\rightarrow}
 \mathcal{N}\left(0,1\right),
$$
which yields, using the second equality in Lemma \ref{Momfin1},
\begin{eqnarray}
\sqrt{nb_1} \left( f_n(\epsilon) - \esp f_{n}(\epsilon) \right)
\stackrel{d}{\rightarrow} \mathcal{N} \left( 0, \frac{f(\epsilon)}
{\prob\left(X\in\mathcal{X}_0\right)} \int K_1^2(v) dv \right).
\label{develop1}
\end{eqnarray}
Moreover, note that for $nb_0^db_1^3\rightarrow\infty$ and
$nb_0^{2d}\rightarrow\infty$,
$$
 \frac{1}{nb_1^5}
 \left(
 \frac{1}{nb_0^d}
 \right)^2
 +
 \left(
 \frac{1}{b_1^2}
 +
 \frac{b_0^d}{b_1^7}
 \right)^2
 \left(\frac{1}{nb_0^d}\right)^3
 =
 O\left(
\frac{1}{n^2b_0^db_1^3}\right).
$$
Therefore, since by Assumptions $(\rm{A}_{11})$ and $(A_9)$, we
have $b_0^4=O\left(1/(nb_0^d)\right)$,
$nb_0^db_1^3\rightarrow\infty$ and that
$nb_0^{2d}\rightarrow\infty$, the equality above and
(\ref{fnchapfn1}) then give
\begin{eqnarray*}
\widehat{f}_{1n}(\epsilon)-\widetilde{f}_{1n}(\epsilon)
 &\asymp&
 O_{\prob}
 \left[
 b_0^4
 +
 \frac{1}{n}
 +
 \frac{1}{n^2b_0^db_1^3}
 +
 \left(
 \frac{1}{nb_1^5}
 +
 \frac{b_0^d}{b_1^3}
 \right)
 \left(
 \frac{1}{nb_0^d}
 \right)^2
 +
 \left(
 \frac{1}{b_1^2}
 +
 \frac{b_0^d}{b_1^7}
 \right)
 \left(\frac{1}{nb_0^d}\right)^3
 \right]^{1/2}
 \\
 &\asymp&
 O_{\prob}
 \left(
 b_0^4
 +
 \frac{1}{n}
 +
 \frac{1}{n^2b_0^db_1^3}
\right)^{1/2}.
\end{eqnarray*}
 Hence for $b_1$ going to $0$, we have
$$
\sqrt{nb_1} \left(
\widehat{f}_{1n}(\epsilon)-\widetilde{f}_{1n}(\epsilon) \right)
 =
O_{\prob}
 \left[
 nb_1
 \left(
 b_0^4
 +
 \frac{1}{n}
 +
 \frac{1}{n^2b_0^db_1^3}
 \right)
\right]^{1/2}
 =
 o_{\prob}(1),
$$
since $nb_0^4b_1=o(1)$ and that $nb_0^db_1^2\rightarrow\infty$
under Assumption $(\rm{A}_{11})$. Combining the above result with
(\ref{develop1}) and (\ref{develop}), we obtain
$$
\sqrt{nb_1}
 \left(
 \widehat{f}_{1n}(\epsilon)
 -
 \esp f_n(\epsilon)
\right) \stackrel{d}{\rightarrow} \mathcal{N}
 \left(
 0,
\frac{f(\epsilon)} {\prob\left(X\in\mathcal{X}_0\right)} \int
K_1^2(v) dv
 \right).
$$
This ends the proof the Theorem, since the first result of Lemma
\ref{Momfin1} gives
$$
\esp f_n(\epsilon) = \esp f_{1n}(\epsilon) = f(\epsilon) +
\frac{b_1^2}{2}f^{(2)}(\epsilon) \int v^2K_1(v) dv +
o\left(b_1^2\right) := \overline{f}_{1n}(\epsilon). \eop
$$

 \setcounter{equation}{0} \setcounter{subsection}{0}
\setcounter{lem}{0}
\renewcommand{\theequation}{A.\arabic{equation}}
\renewcommand{\thesubsection}{A.\arabic{subsection}}
\begin{center}
\section*{Appendix A:  Proof of the  intermediate results}
\end{center}

\subsection*{Proof of Lemma \ref{Estig}}
First note that  by $(A_7)$, we have
 $
 \int\!
z K_0(z) dz
 =0
 $
 and
 $
 \int\!
K_0(z) dz
 =1$.  Therefore, since $K_0$ is continuous and has a compact support,
$(A_1)$, $(A_2)$ and a second-order Taylor expansion,
 yield,
for  $b_0$ small enough and any $x$ in $\mathcal{X}_0$,
\begin{eqnarray*}
 \lefteqn{
 \left|
 \overline{g}_{n}(x)-g(x)
 \right|
  =
 \left|
 \frac{1}{b_0^d}
 \int
 K_0
 \left(\frac{z-x}{b_0}\right)
 g(z)
 dz
 -
 g(x)
 \right|
  =
 \left|
\int K_0(z) \left[g(x+b_0z)-g(x)\right]
 dz
 \right|
 }
 &&
\\
&=&
 \left|
 \int K_0(z)
 \left[
 b_0 g^{(1)}(x) z
 +
 \frac{b_0^2}{2}
 z g^{(2)} (x + \theta b_0 z)z^{\top}
\right] dz
 \right|,
 \;\theta = \theta (x,b_0 z)\in [0,1]
\\
&=& \left|
 b_0 g^{(1)}(x)
 \int
 z K_0(z) dz
 +
 \frac{b_0^2}{2}
  \int
z g^{(2)}(x + \theta b_0 z) z^{\top} K_0(z) dz \right|
\\
& = & \frac{b_0^2}{2}
 \left|
 \int z g^{(2)}(x+\theta b_0z)
 z^{\top} K_0(z) dz
 \right|
 \leq C b_0^2,
\end{eqnarray*}
so that
$$
 \sup_{x\in\mathcal{X}_0}
 \left|
 \overline{g}_n (x) - g(x)
 \right|
 =
 O\left(b_0^2\right).
$$
This gives the first equality of the lemma. To prove the two last
equalities in the Lemma, note that it is
 sufficient to show that
$$
 \sup_{x\in\mathcal{X}_0}
 \left|
 \widehat{g}_{n}(x)
  -
 \overline{g}_n(x)
 \right|
 =
 O_{\prob}
 \left(
  \frac{\ln n}{nb_0^d}
 \right)^{1/2},
$$
 since $\bar{g}_n (x)$ is asymptotically bounded away from
$0$ over
 $\mathcal{X}_0$
 and that
$|\overline{g}_n (x) - g(x) | = O (b_0^2) $ uniformly for $x$ in
$\mathcal{X}_0$. This follows from Theorem 1 in  Einmahl and Mason
(2005). \eop

\subsection*{Proof of Lemma \ref{Estim}}

For the first equality in the lemma, set
$$
\widehat{r}_n(x)
 =
 \frac{1}{nb_0^d}
 \sum_{j=1}^n
 Y_j K_0
  \left(
 \frac{X_j-x}{b_0}
 \right),
 \quad
 \overline{r}_n(x)
  =
  \esp
 \left[
\widehat{r}_n(x)
 \right] \;,
$$
and observe that
\begin{equation}
\sup_{x\in\mathcal{X}_0}
 \left|\widehat{m}_n(x)-m(x)\right|
 \leq
 \sup_{x\in\mathcal{X}_0}
 \left|
 \widehat{m}_n(x)
 -
 \frac{\overline{r}_n(x)}{\overline{g}_n(x)}
 \right|
 +
 \sup_{x\in\mathcal{X}_0}
 \frac{1}{\left|\overline{g}_n(x)\right|}
 \left|
 \overline{r}_n(x)-\overline{g}_n(x)m(x)
 \right|.
 \label{mchd}
\end{equation}
Consider the first term of (\ref{mchd}). Note that
$\esp^{1/4}\left[Y^4|X=x\right] \leq |m(x)|
+\esp^{1/4}\left[\varepsilon^4\right]$. The compactness of
$\mathcal{X}$
 from $(A_1)$, the continuity of $m(\cdot)$ from $(A_3)$
and $(A_4)$ then give that $\esp\left[Y^4|X=x\right]<\infty$
uniformly for $x\in\mathcal{X}_0$. Hence under $(A_9)$, Theorem 2
in Einmahl and Mason (2005) gives
$$
\sup_{x\in\mathcal{X}_0} \left|
 \widehat{m}_n(x)
 -
\frac{\overline{r}_n(x)}{\overline{g}_n(x)}
 \right|
 =
 O_{\prob}
 \left(\frac{\ln n}{nb_0^d}\right)^{1/2}.
 $$
 For the second term in (\ref{mchd}), a
second-order Taylor expansion
 gives, as in the proof of Lemma \ref{Estig},
$$
 \sup_{x\in\mathcal{X}_0}
 \left|
 \overline{r}_n(x)
 -
 \overline{g}_n(x)m(x)
\right|
 =
 O (b_0^2).
 $$
  This gives the result of lemma since
 Lemma \ref{Estig}
 implies that $\overline{g}_n (x)$ is bounded away from $0$ over
 $\mathcal{X}_0$ uniformly in $x$ and for $b_0$ small enough.
 \eop

\subsection*{Proof of Lemma \ref{Reste}}
Note that under $(A_8)$, the Taylor expansion with integral
remainder gives, for any $x\in\mathcal{X}_0$ and any integer
$i\in[1,n]$,
\begin{eqnarray*}
 K_1\left(
 \frac{Y_i-\widehat{m}_n(x)-\epsilon}{h_1}
 \right)
=
 K_1\left(\frac{Y_i-m(x)-\epsilon}{h_1}\right)
 -
 \frac{1}{h_1}
 \left(\widehat{m}_n(x)-m(x)\right)
 \int_{0}^1
 K_1^{(1)}
 \left(\frac{Y_i-\theta_n(x, t)}{h_1}
 \right)
  dt,
\end{eqnarray*}
where
 $\theta_n(x,t)=m(x)+\epsilon+t\left(\widehat{m}_n(x)-m(x)\right)$.
 Therefore
 \begin{eqnarray}
 \nonumber
 \widetilde{f}_n(\epsilon|x)
 =
 f_{n}(\epsilon|x)
 -
 \frac{\widehat{m}_n(x)-m(x)}{\widetilde{g}_n(x)}
 \left[
 \frac{1}{nh_0^dh_1^2}
 \sum_{i=1}^n
 K_0\left(\frac{X_i-x}{h_0}\right)
 \int_0^1
 K_1^{(1)}
 \left(\frac{Y_i-\theta_n(x,t)}{h_1}\right)
 dt
 \right].
 \\
 \label{fnt}
 \end{eqnarray}
Now, observe that if $X_i = z$ and $y\in\Rit$, the change of
variable $e=y-m(z)+h_1v$ gives,  under  $(A_1)-(A_5)$ and $(A_7)$,
\begin{eqnarray*}
\lefteqn{ \esp_n
 \left|
 K_1^{(1)}
 \left(\frac{Y_i-y}{h_1}\right)
 \right|
 =
 \esp
 \left|
 K_1^{(1)}
 \left(\frac{\varepsilon_i+m(z)-y}{h_1}\right)
\right|
 }
 \\
 &=&
 \int
 \left|
 K_1^{(1)}
 \left(\frac{e+m(z)-y}{h_1}\right)
\right|
 f(e) de
 \\
 &=&
 h_1
 \int
|K_1^{(1)}(v)|
 f\left((y-m(z)+h_1v\right))
 dv
\leq Ch_1.
\end{eqnarray*}
Hence
$$
\sup_{1\leq i\leq n}
 \int_0^1
 \esp_n
 \left|
 K_1^{(1)}
 \left(\frac{Y_i-\theta_n(x,t)}{h_1}\right)
 \right|
 dt
\leq Ch_1.
$$
With the help of this result and Lemma \ref{Estig}, we have
\begin{eqnarray*}
\lefteqn{ \esp_n \left| \frac{1}{nh_0^dh_1}
 \sum_{i=1}^n
 K_0\left(\frac{X_i-x}{h_0}\right)
 \int_0^1
 K_1^{(1)}
 \left(\frac{Y_i-\theta_n(x,t)}{h_1}\right)
 dt
\right| }
\\
&\leq& \frac{1}{nh_0^dh_1}
 \sum_{i=1}^n
 \left|K_0\left(\frac{X_i-x}{h_0}\right)\right|
 \times
\sup_{1\leq i\leq n} \int_0^1
 \esp_n
 \left|
 K_1^{(1)}
 \left(\frac{Y_i-\theta_n(x,t)}{h_1}\right)
 \right|
 dt
 \\
 &\leq&
\frac{C}{nh_0^d}
 \sum_{i=1}^n
 \left|
 K_0\left(\frac{X_i-x}{h_0}\right)
 \right|
 =
 O_{\prob}(1),
\end{eqnarray*}
so that
$$
\frac{1}{nh_0^dh_1^2}
 \sum_{i=1}^n
 K_0\left(\frac{X_i-x}{h_0}\right)
 \int_0^1
 K_1^{(1)}
 \left(\frac{Y_i-\theta_n(x,t)}{h_1}\right)
 dt
=O_{\prob}\left(\frac{1}{h_1}\right).
$$
Hence from (\ref{fnt}),  (\ref{gtilde}),  Lemma \ref{Estim} and
Assumption $(\rm\bf{A}_0)$, we deduce
$$
 \widetilde{f}_n(\epsilon|x)
 =
 f_{n}(\epsilon|x)
 +
 O_{\prob}
 \left(\frac{1}{h_1}\right)
 \left(b_0^4+\frac{\ln n}{nb_0^d}\right)^{1/2}
 =
 f_{n}(\epsilon|x)
 +
 o\left(\frac{1}{nh_0^dh_1}\right)^{1/2}.
 \eop
$$

\subsection*{Proof of Lemma \ref{Momvarphi} and Lemma \ref{Momfin1}}
We just give the proof of Lemma \ref{Momvarphi}, the proof of
Lemma \ref{Momfin1} being very similar.
 For the first equality of  Lemma \ref{Momvarphi}, note that
\begin{eqnarray*}
\esp \left[ \widetilde{\varphi}_{in}(x,y) \right] &=& \nonumber
\frac{1} {h_0^dh_1} \int \int K_0 \left( \frac{x_1-x}{h_0} \right)
K_1 \left( \frac{y_1-y}{h_1} \right) \varphi (x_1, y_1) dx_1 dy_1
\\
&=& \int \int K_0(z) K_1(v) \varphi\left(x+h_0z, y+h_1 v\right) dz
dv.
 \label{espphi}
\end{eqnarray*}
A second-order Taylor expansion gives under $(A_6)$, for $z$ in
the support of $K_0$, $v$ in the support of $K_1$, and $h_0$,
$h_1$ small enough,
\begin{eqnarray*}
\lefteqn{ \varphi\left(x+h_0z, y+h_1 v\right) - \varphi (x, y) }
\\
&=& h_0 \frac{\partial \varphi (x, y)}{\partial x} z^{\top} + h_1
\frac{\partial \varphi (x, y)}{\partial y} v
\\
&& + \frac{h_0^2}{2} z \frac{\partial^2 \varphi (x+\theta h_0 z,
y+\theta h_1 v)}{\partial^2
 x} z^{\top}
+ h_1 h_0 v \frac{\partial^2 \varphi (x+\theta h_0 z, y+\theta h_1
v)}{\partial x\partial y} z^{\top}
\\
&& + \frac{h_1^2}{2} \frac{\partial^2 \varphi (x+\theta h_0 z,
y+\theta h_1 v)}{\partial^2
 y}
v^2,
\end{eqnarray*}
for some $\theta = \theta (x,y,h_0 z, h_1 v)$ in $[0,1]$. This
gives, since $\int\! K_0 (z) dz = \int\! K_1 (v) dv =1$, $\int\!z
K_0(z) dz$ and $\int\! v K_1 (v) dv$ vanish under $(A_7)-(A_8)$,
and by the Lebesgue Dominated Convergence Theorem,
\begin{eqnarray*}
\lefteqn{ \esp \left[ \widetilde{\varphi}_{in}(x,y) \right] -
\varphi(x, y) - \frac{h_0^2}{2} \frac{\partial^2 \varphi(x,
y)}{\partial^2 x} \int z K_0(z)z^{\top} dz - \frac{h_1^2}{2}
\frac{\partial^2 \varphi (x, y)}{\partial^2 y} \int v^2 K_1 (v) dv
} &&
\\
& = & \frac{h_0^2}{2} \int \int z \left( \frac{\partial^2 \varphi
(x+\theta h_0 z, y+\theta h_1 v)}{\partial^2
 x}
- \frac{\partial^2 \varphi (x, y)}{\partial^2 x} \right) z^{\top}
K_0 (z) K_1 (v) dz dv
\\
&& + h_1 h_0 \int \int v \left( \frac{\partial^2 \varphi (x+\theta
h_0 z, y+\theta h_1 v)} {\partial x \partial y} - \frac{\partial^2
\varphi (x, y)}{\partial x \partial y} \right) z^{\top} K_0 (z)
K_1 (v) dz dv
\\
&& + \frac{h_1^2}{2} \int \int \left( \frac{\partial^2 \varphi
(x+\theta h_0 z, y+\theta h_1 v)} {\partial^2 y} -
\frac{\partial^2 \varphi (x, y)}{\partial^2 y} \right) v^2 K_0 (z)
K_1 (v) dz dv
\\
&=& o ( h_0^2 + h_1^2).
\end{eqnarray*}
This proves the first equality of Lemma \ref{Momvarphi}. The
second equality in Lemma  follows similarly, since
\begin{eqnarray*}
\lefteqn{ \Var [\widetilde{\varphi}_{in}(x, y)] = \esp \left[
\widetilde{\varphi}_{in}^2(x, y) \right] - \left( \esp \left[
\widetilde{\varphi}_{in}(x, y) \right] \right)^2 }
\\
& = & \frac{1}{h_0^dh_1} \int \int \varphi\left(x+h_0z,
y+h_1v\right) K_0^2(z) K_1^2(v) dz dv + O(1)
\\
&=& \frac{\varphi(x, y)} {h_0^dh_1} \int \int K_0^2(z) K_1^2(v) dz
dv + o\left(\frac{1}{h_0^dh_1}\right).
\end{eqnarray*}

The last statement of Lemma \ref{Momvarphi} is immediate, since
the Triangular and Convex inequalities give
\begin{eqnarray*}
\esp \left| \widetilde{\varphi}_{in}(x, y) - \esp
\widetilde{\varphi}_{in}(x, y) \right|^3 &\leq& C\esp \left|
\widetilde{\varphi}_{in}(x, y) \right|^3
\\
&\leq& \frac {C\varphi(x, y)} {h_0^{2d}h_1^2 } \int \int \left|
K_0(z) K_1(v) \right|^3 dz dv +
o\left(\frac{1}{h_0^{2d}h_1^2}\right). \eop
\end{eqnarray*}

\subsection*{Proof of Lemma \ref{STR}}
The order of $S_n$ follows from Lemma \ref{Betasum} and Lemma
\ref{Sigsum}. In fact,  since
\begin{eqnarray*}
\mathds{1}(X_i \in \mathcal{X}_0)
 \left(
 \widehat{m}_{in} - m(X_i)
\right)
 &=&
 \frac{\mathds{1}(X_i \in \mathcal{X}_0)}
 {nb_0^d\widehat{g}_{in}}
 \sum_{j=1,j\neq i}^n
 \left( m(X_j)+ \varepsilon_j - m(X_i)\right)
 K_0\left(\frac{X_j-X_i}{b_0}\right)
 \\
 & =&
 \beta_{in} + \Sigma_{in},
\end{eqnarray*}
Lemma \ref{Betasum} and Lemma \ref{Sigsum} give
\begin{eqnarray*}
S_n
 =
 O_{\prob}
\left[ b_0^2 \left(nb_1^2+(nb_1)^{1/2}\right)
 +
\left( nb_1^4 + \frac{b_1}{b_0^d} \right)^{1/2}
 \right],
\end{eqnarray*}
which gives the result for $S_n$.

\vskip 0.3cm For $T_n$, define for any $1\leq i\leq n$,
$$
\esp_{in} [\cdot] = \esp_n \left[
 X_1,\ldots,X_n,\varepsilon_k,
 k\neq i
 \right].
$$
Therefore, since $(\widehat{m}_{in} - m(X_i))$ depends only upon
$\left(X_1,\ldots,X_n,\varepsilon_k,k\neq i\right)$, we have
\begin{eqnarray*}
 \esp_n [T_{n}]
 &=&
 \esp_{n}
 \left[
 \sum_{i=1}^n
 \esp_{in}
 \left[
 \mathds{1}
\left( X_i \in \mathcal{X}_0 \right)
 (\widehat{m}_{in} - m(X_i))^2
  K_1^{(2)}
 \left(\frac{\varepsilon_{i}-\epsilon}{b_1}\right)
 \right]
 \right]
 \\
 &=&
 \esp_{n}
  \left[
  \sum_{i=1}^n
  \mathds{1}
 \left( X_i \in \mathcal{X}_0\right)
(\widehat{m}_{in} - m(X_i))^2
 \esp_{in}
 \left[
 K_1^{(2)}
 \left(\frac{\varepsilon_{i}-\epsilon}{b_1}\right)
 \right]
 \right],
\end{eqnarray*}
with, using $(A_4)$ and Lemma \ref{MomderK}-(\ref{MomderK2}),
\begin{eqnarray*}
 \left|
\esp_{in}
 \left[
 K_1^{(2)}
 \left(\frac{\varepsilon_{i}-\epsilon}{b_1}\right)
 \right]
 \right|
 =
 \left|
 \int
  K_1^{(2)}
 \left(\frac{e-\epsilon}{b_1}
 \right)
 f(e)
 de
 \right|
\leq
 Cb_1^3.
\end{eqnarray*}
 Hence this bound, the equality above, the Cauchy-Schwarz inequality
  and  Lemma \ref{BoundEspmchap} yield that
\begin{eqnarray}
\nonumber
 \left|
 \esp_{n}\left[T_{n}\right]
 \right|
&\leq&
 Cb_1^3
 \sum_{i=1}^n
  \esp_{n}
  \biggl[
  \mathds{1}\left( X_i \in \mathcal{X}_0\right)
 (\widehat{m}_{in} - m(X_i))^2
 \biggr]
 \\\nonumber
 &\leq&
 C n b_1^3
 \left(
 \sup_{1\leq i\leq n}
 \esp_{n}
  \biggl[
  \mathds{1}\left( X_i \in \mathcal{X}_0\right)
 (\widehat{m}_{in} - m(X_i))^4
 \biggr]
 \right)^{1/2}
 \\
 &\leq&
 O_{\prob}
 \left(nb_1^3\right)
 \left(b_0^4+\frac{1}{nb_0^d}\right).
 \label{Boundmean}
\end{eqnarray}

For the conditional variance of $T_n$,  Lemma \ref{sumzeta} gives
\begin{eqnarray*}
\Var_n (T_{n})
 &=&
 \sum_{i=1}^n
 \Var_n\left( \zeta_{in}\right)
 +
 \sum_{i=1}^n
 \sum_{j=1\atop j\neq i}^n
\Cov_n \left( \zeta_{in} , \zeta_{jn}
 \right)
\\
&=& O_{\prob} \left(nb_1\right)
\left(b_0^4+\frac{b_1}{nb_0^d}\right)^2 + O_{\prob}
\left(n^2b_0^db_1^{7/2}\right)
\left(b_0^4+\frac{1}{nb_0^d}\right)^2.
\end{eqnarray*}
Therefore, since $b_1$ goes to $0$ under $(A_{10})$,  this order
and (\ref{Boundmean}) yield,
 applying the  Tchebychev inequality,
\begin{eqnarray*}
T_{n}
 &=&
O_{\prob} \left[ \left(nb_1^3\right)
\left(b_0^4+\frac{1}{nb_0^d}\right)
 +
\left(nb_1\right)^{1/2} \left(b_0^4+\frac{b_1}{nb_0^d}\right) +
\left(n^2b_0^db_1^{7/2}\right)^{1/2}
\left(b_0^4+\frac{1}{nb_0^d}\right) \right]
\\
&=& O_{\prob} \left[ \left( nb_1^3 + \left(nb_1\right)^{1/2} +
\left(n^2b_0^db_1^{3}\right)^{1/2} \right)
\left(b_0^4+\frac{1}{nb_0^d}\right) \right].
\end{eqnarray*}
which gives the result for $T_n$.

\vskip 0.2cm We now compute the order of $R_n$. For this, define
\begin{eqnarray*}
I_{in} &= & \int_{0}^{1}
 (1-t)^2
 K_1^{(3)}
 \left(
 \frac{
 \varepsilon_i-t(\widehat{m}_{in} -m(X_i))-\epsilon}{b_1}
 \right)
 dt,
 \\
R_{in} & =&
 \mathds{1}
 \left(X_i\in\mathcal{X}_0\right)
 \left(\widehat{m}_{in}-m(X_i)\right)^3
 I_{in},
\end{eqnarray*}
and note that $R_n=\sum_{i=1}^n R_{in}$. The order of $R_n$ is
derived by computing its conditional mean and its conditional
variance. For the conditional mean,  observe that
\begin{eqnarray*}
 \esp_n [R_{n}]
&=&
 \esp_n
 \left[
 \sum_{i=1}^n
 \esp_{in}\left[R_{in}\right]
\right]
\\
 &=&
 \esp_{n}
 \left[
 \sum_{i=1}^n
 \mathds{1}
 \left( X_i \in \mathcal{X}_0\right)
 (\widehat{m}_{in} - m(X_i))^3
 \esp_{in}
 \left[ I_{in}\right]\right],
\end{eqnarray*}
with, using $(A_4)$ and Lemma \ref{MomderK}-(\ref{MomderK3}),
\begin{eqnarray*}
 \left|
 \esp_{in}
 \left[ I_{in}\right]
 \right|
& =&
 \left|
 \int_{0}^{1}
 (1-t)^2
 \left[
 \int
  K_1^{(3)}
 \left(
 \frac{e-t(\widehat{m}_{in} -m(X_i))-\epsilon}{b_1}
 \right)f(e)
 de
 \right]
 dt
 \right|
 \\
&\leq&
 Cb_1^3.
\end{eqnarray*}
Therefore the Holder inequality and Lemma \ref{BoundEspmchap}
yield
\begin{eqnarray}
\nonumber
 \left|
 \esp_n\left[R_n\right]
 \right|
 &\leq&
 Cb_1^3
 \sum_{i=1}^n
 \esp_n
 \left[
\left| \mathds{1} \left(X_i\in\mathcal{X}_0\right)
\left(\widehat{m}_{in}-m(X_i)\right) \right|^3 \right]
\\\nonumber
&\leq& Cb_1^3 \sum_{i=1}^n \esp_n^{3/4} \left[ \mathds{1}
\left(X_i\in\mathcal{X}_0\right)
\left(\widehat{m}_{in}-m(X_i)\right)^4 \right]
\\
&\leq&
 O_{\prob}
 \left(nb_1^3\right)
 \left(b_0^4+\frac{1}{nb_0^d}\right)^{3/2}.
 \label{Rn1}
\end{eqnarray}
For the conditional covariance of $R_n$, note that Lemma
\ref{Indep} allows to write
\begin{eqnarray}
\Var_n\left(R_n\right) = \sum_{i=1}^n \Var_n\left(R_{in}\right) +
\sum_{i=1}^n \sum_{j=1\atop j\neq i}^n \biggl(\|X_i-X_j\|\leq
Cb_0\biggr) \Cov_n\left(R_{in}, R_{jn}\right), \label{Rn2}
\end{eqnarray}
and consider  the  first term in (\ref{Rn2}). We have
\begin{eqnarray*}
\Var_n\left(R_{in}\right)
 \leq
\esp_n\left[R_{in}^2\right] \leq
 \esp_{n}
 \biggl[
 \mathds{1}
 \left( X_i \in \mathcal{X}_0\right)
 (\widehat{m}_{in} - m(X_i))^6
 \esp_{in}
 \left[ I_{in}^2\right]
 \biggr],
\end{eqnarray*}
with, using  $(A_4)$, the Cauchy-Schwarz inequality and Lemma
\ref{MomderK}-(\ref{MomderK3}),
\begin{eqnarray*}
 \esp_{in}
 \left[ I_{in}^2\right]
 &\leq&
 C\esp_{in}
 \left[
\int_{0}^{1}
 K_1^{(3)}
 \left(
 \frac{\varepsilon_i-t(\widehat{m}_{in} -m(X_i))-\epsilon}{b_1}
 \right)^2
 dt
 \right]
\\
&\leq&
 C\int_{0}^{1}
 \left[
 \int
 K_1^{(3)}
 \left(
 \frac{e-t(\widehat{m}_{in} -m(X_i))-\epsilon}{b_1}
 \right)^2
 f(e)
 de
 \right]
 dt
 \\
&\leq&
 Cb_1,
\end{eqnarray*}
so that
\begin{eqnarray*}
 \Var_n\left(R_{in}\right)
\leq
 Cb_1
 \esp_{n}
 \left[
 \mathds{1}
 \left( X_i \in \mathcal{X}_0\right)
 (\widehat{m}_{in} - m(X_i))^6
 \right].
 \end{eqnarray*}
Therefore form Lemma \ref{BoundEspmchap}, we deduce
\begin{eqnarray}
\nonumber
 \sum_{i=1}^n
 \Var_n\left(R_{in}\right)
 &\leq&
 Cnb_1
 \sup_{1\leq i\leq n}
 \esp_{n}
 \left[
 \mathds{1}
 \left( X_i \in \mathcal{X}_0\right)
 (\widehat{m}_{in} - m(X_i))^6
 \right]
 \\
&\leq&
 O_{\prob}
 \left(nb_1\right)
 \left(b_0^4+\frac{1}{nb_0^d}\right)^3.
 \label{Rn4}
\end{eqnarray}
For the second  term in (\ref{Rn2}),  the Cauchy-Schwarz
inequality gives, with the help of the above result for
$\Var_n\left(R_{in}\right)$,
\begin{eqnarray*}
\left| \Cov_n\left(R_{in}, R_{jn}\right)
 \right|
 &\leq&
 \left(
 \Var_n\left(R_{in}\right)
 \Var_n\left(R_{jn}\right)
 \right)^{1/2}
 \\
 &\leq&
 Cb_1
 \sup_{1\leq i\leq n}
 \esp_{n}
 \left[
 \mathds{1}
 \left( X_i \in \mathcal{X}_0\right)
 (\widehat{m}_{in} - m(X_i))^6
 \right].
\end{eqnarray*}
Hence by Lemma \ref{BoundEspmchap}  and the Markov inequality, we
have
\begin{eqnarray*}
\lefteqn{ \sum_{i=1}^n \sum_{j=1\atop j\neq i}^n
\biggl(\|X_i-X_j\|\leq Cb_0\biggr) \left| \Cov_n\left(R_{in},
R_{jn}\right) \right| }
\\
&\leq&
 O_{\prob}\left(b_1\right)
 \left(b_0^4+\frac{1}{nb_0^d}\right)^3
 \sum_{i=1}^n
 \sum_{j=1\atop j\neq i}^n
 \biggl(\|X_i-X_j\|\leq Cb_0\biggr)
 \\
 &\leq&
 O_{\prob}\left(b_1\right)
 \left(b_0^4+\frac{1}{nb_0^d}\right)^3
 \left(n^2b_0^d\right).
\end{eqnarray*}
This order, (\ref{Rn4}) and (\ref{Rn2}) give, since $nb_0^d$
diverges  under $(A_9)$,
$$
\Var\left(R_n\right) =
 O_{\prob}
 \left(b_0^4+\frac{1}{nb_0^d}\right)^3
 \left(n^2b_0^db_1\right).
$$
Finally, with the help of this result, (\ref{Rn1}) and the
Tchebychev inequality, we arrive at
\begin{eqnarray*}
R_n &=& O_{\prob} \left[ \left(nb_1^3\right)
 \left(b_0^4+\frac{1}{nb_0^d}\right)^{3/2}
 +
\left(n^2b_0^db_1\right)^{1/2}
\left(b_0^4+\frac{1}{nb_0^d}\right)^{3/2} \right]
\\
&=& O_{\prob} \left[ \left( nb_1^3 +
\left(n^2b_0^db_1\right)^{1/2} \right)
\left(b_0^4+\frac{1}{nb_0^d}\right)^{3/2} \right].
 \eop
\end{eqnarray*}

\subsection*{Proof of Lemma \ref{MomderK}}
Set $h_p(e)=e^p f(e)$, $p\in[0,2]$. For  the first inequality of
(\ref{MomderK1}), note that  under $(A_5)$ and $(A_8)$, the change
of variable $e=\epsilon+b_1 v$ give, for any integer $\ell\in[1,
3]$,
\begin{eqnarray}
\nonumber
 \left|
 \int
 K_1^{(\ell)}
 \left(\frac{e-\epsilon}{b_1}\right)^2
 e^pf(e) de
 \right|
 &=&
 \left|
 b_1
 \int
 K_1^{(\ell)}(v)^2 h_p(\epsilon+b_1v)
 dv
 \right|
 \\\nonumber
 &\leq&
 b_1
 \sup_{t\in\Rit}
 |h_p(t)|
 \int
 | K_1^{(\ell)}(v)^2|
 dv
 \\
 &\leq&
 Cb_1,
 \label{Ineg1}
\end{eqnarray}
which yields the first inequality in (\ref{MomderK1}). For the
second inequality in (\ref{MomderK1}), observe that $f(\cdot)$ has
a bounded  continuous derivative under $(A_5)$, and that $\int \!
K_1^{(\ell)}(v)dv =0$ under $(A_8)$. Therefore, since $h_p(\cdot)$
has bounded second order derivatives under $(A_7)$, the Taylor
inequality yields that
\begin{eqnarray*}
\left| \int
 K_1^{(\ell)}
\left(\frac{e-\epsilon}{b_1}\right)
 e^pf(e)de
 \right|
 &=&
  b_1
 \left|
 \int
 K_1^{(\ell)}(v)
 \left[
 h_p(\epsilon+b_1v)-h_p(\epsilon)
 \right]
 \right|
  dv
 \\
 &\leq&
 b_1^2
 \sup_{t\in\Rit}|h_p^{(1)}(t)|
 \int
 |vK_1^{(\ell)}(v)|
  dv
 \leq
 Cb_1^2.
\end{eqnarray*}
which  completes the proof of (\ref{MomderK1}).

 The first inequalities of  (\ref{MomderK2}) and
 (\ref{MomderK3}) follow directly from (\ref{Ineg1}). The second bounds in
(\ref{MomderK2}) and (\ref{MomderK3}) are proved simultaneously.
For this, note that for any integer $\ell\in\{2,3\}$,
$$
\int
 K_1^{(\ell)}
 \left(\frac{e-\epsilon}{b_1}\right)
 h_p(e) de
 =
b_1 \int
 K_1^{(\ell)}(v)
h_p(\epsilon+b_1v) dv.
$$
Under $(A_8)$,  $K_1(\cdot)$ is symmetric, has  a compact support
and two
 continuous derivatives, with
$\int \! K_1^{(\ell)}(v)dv=0$ and
 $\int\! v K_1^{(\ell)}(v) dv=0$ for $\ell\in\{2,3\}$.
 Hence, since by $(A_5)$ $h_p$ has bounded continuous second
order derivatives, this gives for some $\theta=\theta
(\epsilon,b_1 v)$,
\begin{eqnarray*}
 \lefteqn{
 \left|
 \int
 K_1^{(\ell)}
\left(
 \frac{e-\epsilon}{b_1}
 \right)
  h_p(e) de
 \right|
 =
 \left|
  b_1
 \int
 K_1^{(\ell)}(v)
\left[
 h_p(\epsilon+b_1v) - h_p(\epsilon)
 \right]
 dv
 \right|
 }
 &&
\\
&=& \left| b_1 \int
 K_1^{(\ell)}(v)
\left[
 b_1 v h_p^{(1)}(\epsilon)
 +
  \frac{b_1^2v^2}{2}
h_p^{(2)}(\epsilon+\theta b_1v) \right]
 dv
 \right|
 \\
  &=&
 \left|
\frac{b_1^3}{2} \int v^2 K_1^{(\ell)}(v)
 h_p^{(2)}(\epsilon+\theta b_1v)
 dv
\right|
\\
&\leq&
 \frac{b_1^3}{2}
 \sup_{t\in \Rit}|h_p^{(2)}(t)|
 \int
\left| v^2K_1^{(\ell)}(v)
 \right|
 dv
 \leq
 Cb_1^3.
 \eop
\end{eqnarray*}

\subsection*{Proof of Lemma \ref{Betasum}}
Assumption $(A_4)$ and Lemma \ref{MomderK}-(\ref{MomderK1}) give
\begin{eqnarray*}
\left| \esp_n \left[
 \sum_{i=1}^n
 \beta_{in}
 K_1^{(1)}
 \left(
\frac{\varepsilon_i-\epsilon}{b_1} \right)
 \right]
 \right|
 & = &
\left| \esp \left[
 K_1^{(1)}
  \left(
  \frac{\varepsilon -\epsilon}{b_1}
  \right)
  \right]
  \sum_{i=1}^n
  \beta_{in}
   \right|
  \leq C n b_1^2
 \max_{1 \leq i \leq n}
 \left| \beta_{in} \right|,
\\
\Var_n
 \left[
 \sum_{i=1}^n
  \beta_{in}
  K_1^{(1)}
  \left(
\frac{\varepsilon_i-\epsilon}{b_1}
 \right)
 \right]
&\leq&
  \sum_{i=1}^n
 \beta_{in}^2
 \esp
 \left[
 K_1^{(1)}
  \left(
   \frac{\varepsilon - \epsilon}{b_1}
\right)^2
 \right]
  \leq
 C n b_1
  \max_{1 \leq i \leq n}
  \left|
  \beta_{in}
\right|^2 .
\end{eqnarray*}
Hence the (conditional) Markov inequality gives
$$
\sum_{i=1}^n \beta_{in}
 K_1^{(1)}
 \left(
\frac{\varepsilon_i-\epsilon}{b_1} \right) = O_{\prob}
 \left( n
b_1^2 + (nb_1)^{1/2}
 \right)
 \max_{1 \leq i \leq n}
  \left|
\beta_{in}
 \right|,
$$
so that the lemma follows if we can prove that
\begin{equation}
\sup_{1 \leq i \leq n}
 \left|\beta_{in}\right|
 =
 O_{\prob}
\left( b_0^2\right),
 \label{BetasumTBP}
\end{equation}
as established now. For this, define
$$
\zeta_j (x)
 =
 \mathds{1}
 \left( x \in \mathcal{X}_0\right)
\left(m(X_j)-m(x)\right)
 K_0\left( \frac{X_j-x}{b_0}\right),
 \;\;
\nu_{in}(x)
 =
 \frac{1}{(n-1) b_0^d}
 \sum_{j=1, j\neq i}^n
 \left(
\zeta_j (x)-\esp[\zeta_j (x)]
 \right),
$$
and $\bar{\nu}_{n} (x)= \esp[\zeta_j(x)] / b_0^d $, so that
$$
\beta_{in} =
 \frac{n-1}{n} \frac{\nu_{in} (X_i)
 +
 \bar{\nu}_n (X_i)}{\widehat{g}_{in}}
\;.
$$
For  $\max_{1 \leq i \leq n} | \bar{\nu}_n (X_i) |$, first observe
that  a second-order Taylor expansion applied successively to
$g(\cdot)$ and $m(\cdot)$ give, for $b_0$ small enough, and for
any $x$, $z$ in $\mathcal{X}$,
\begin{eqnarray*}
 \lefteqn
 {
 \left[ m(x+b_0z)-m(x)\right]
 g(x+b_0z)
 }
\\
& =& \left[
 b_0 m^{(1)}(x) z
 +
 \frac{b_0^2}{2}
z m^{(2)}(x +\zeta_1 b_0 z)z^{\top}
 \right]
 \left[
 g(x)
 +
 b_0 g^{(1)}(x) z
 +
 \frac{b_0^2}{2}
 z g^{(2)}(x +\zeta_2 b_0z)z^{\top}
 \right],
  \end{eqnarray*}
for some $\zeta_1 =\zeta_1 (x,b_0 z)$ and $\zeta_2 =\zeta_2 (x,b_0
z)$ in $[0,1]$. Therefore, since $\int\!z
 K(z)dz=0$ under $(A_7)$, it follows that, by $(A_1)$, $(A_2)$ and
 $(A_3)$,
\begin{eqnarray}
\nonumber
 \max_{1 \leq i \leq n}
|\bar{\nu}_n (X_i)| &\leq& \sup_{x\in\mathcal{X}_0} |\bar{\nu}_n
(x)|
 =
 \sup_{x \in
 \mathcal{X}_0}
 \left|
 \int
 \left(m ( x + b_0 z) - m(x)\right)
 K_0 (z) g(x+b_0z)
 dz
\right|
\\
&\leq&
 Cb_0^2.
 \label{Betasum1}
\end{eqnarray}
 Consider now the term  $\max_{1\leq i \leq n}|\nu_{in}(X_i)|$.
  The Bernstein inequality (see e.g. Serfling (2002)) and
$(A_4)$ give, for any $t>0$,
\begin{eqnarray*}
\prob \left( \max_{1 \leq i \leq n} | \nu_{in} (X_i)|
 \geq t
\right) &\leq & \sum_{i=1}^n \prob
 \left(
 | \nu_{in} (X_i) |
 \geq
t \right)
 \leq
 \sum_{i=1}^n
 \int
 \prob
 \left(
  | \nu_{in} (x) |
\geq t
 \left| X_i = x \right.
 \right) g (x)
  dx
\\
& \leq &
 2n \exp
 \left(
  - \frac{ (n-1) t^2 }
  { 2\sup_{x \in\mathcal{X}_0}
\Var (\zeta_j (x)/b_0^d) + \frac{4M}{3b_0^d} t}
 \right),
\end{eqnarray*}
where $M$ is such that $\sup_{x \in \mathcal{X}_0}|\zeta_j (x)|
\leq M$. The definition of $\mathcal{X}_0$ given in $(A_2)$,
$(A_3)$, $(A_7)$ and the standard Taylor expansion yield, for
$b_0$ small enough,
$$
\sup_{x \in \mathcal{X}_0}
 | \zeta_j (x) |
 \leq C b_0,
 \;\;\;
\sup_{x \in \mathcal{X}_0}
 \Var (\zeta_j (x)/b_0^d)
\leq \frac{1}{b_0^d}
 \sup_{x \in \mathcal{X}_0}
 \int
 \left( m(x +b_0 z) - m(x) \right)^2
 K_0^2 (z) g(x+b_0z) dz
  \leq
  \frac{C
b_0^2}{b_0^d}\;,
$$
so that, for any $t \geq 0$,
$$
\prob
 \left(
 \max_{1 \leq i \leq n}
 | \nu_{in} (X_i) |
 \geq t
\right) \leq 2n
 \exp
 \left( - \frac{(n-1) b_0^d t^2 /b_0^2}{C + C
t/b_0} \right).
$$
This gives
$$
\prob \left( \max_{1 \leq i \leq n} |\nu_{in} (X_i)| \geq \left(
\frac{b_0^2 \ln n}{ (n-1) b_0^d} \right)^{1/2} t\right) \leq 2n
\exp \left(
 - \frac{ t^2 \ln n }
 {C + C t \left( \frac{\ln n}{
(n-1) b_0^d} \right)^{1/2} } \right) = o(1),
$$
provided that $t$ is large enough and under $(A_9)$. It then
follows that
$$
\max_{1 \leq i \leq n} | \nu_{in} (X_i) |
 =
 O_{\prob}
 \left(
 \frac{b_0^2 \ln n}{ n b_0^d} \right)^{1/2}.
$$
This bound, (\ref{Betasum1}) and Lemma \ref{Estig}  show that
(\ref{BetasumTBP}) is proved,  since $b_0^2\ln
n/(nb_0^d)=O\left(b_0^4\right)$ under $(A_9)$, and that
$$
 \beta_{in}
  =
 \frac{n-1}{n} \frac{\nu_{in} (X_i)
 +
 \bar{\nu}_n (X_i)}{\widehat{g}_{in}}\;.
 \eop
$$

\subsection*{Proof of Lemma \ref{Sigsum}}
Note that $(A_4)$ gives that  $\Sigma_{in}$ is independent of
$\varepsilon_i$, and that $\esp_n[\Sigma_{in}]=0$. This yields
\begin{eqnarray}
 \esp_n
 \left[
 \sum_{i=1}^n
 \Sigma_{in}
 K_1^{(1)}
 \left(
 \frac{\varepsilon_i-\epsilon}{b_1}
 \right)
 \right]
  = 0.
\label{EspSigmai}
\end{eqnarray}
Moreover,  observe that
\begin{eqnarray}
\nonumber
 \lefteqn{
 \Var_n
  \left[
 \sum_{i=1}^n
 \Sigma_{in}
 K_1^{(1)}
  \left(
\frac{\varepsilon_i-\epsilon}{b_1}
 \right)
 \right]
 }
\\\nonumber
&=&
 \sum_{i=1}^n
 \Var_n
  \left[
 \Sigma_{in}
 K_1^{(1)}
 \left(
 \frac{\varepsilon_i-\epsilon}{b_1}
 \right)
 \right]
 +
 \sum_{i=1}^n
 \sum_{j=1\atop j\neq i}^n
 \Cov_n
 \left[
 \Sigma_{in}
 K_1^{(1)}
 \left(
 \frac{\varepsilon_{i}-\epsilon}{b_1}
 \right)
 ,
 \Sigma_{jn}
 K_1^{(1)}
 \left(
 \frac{\varepsilon_{j}-\epsilon}{b_1}
 \right)
 \right].
 \\
\label{VarSigm}
\end{eqnarray}
For the sum of variances in (\ref{VarSigm}), Lemma
\ref{MomderK}-(\ref{MomderK1}) and $(A_4)$ give
\begin{eqnarray}
\nonumber \sum_{i=1}^n
 \Var_n
 \left[
 \Sigma_{in} K_1^{(1)}
 \left(
 \frac{\varepsilon_i-\epsilon}{b_1}
 \right)
 \right]
 &\leq&
 \sum_{i=1}^n
 \esp_n
\left[ \Sigma_{in}^2
 \right]
 \esp
 \left[
 K_1^{(1)}
  \left(
\frac{\varepsilon_i-\epsilon}{b_1} \right)^2
 \right]
\\\nonumber
&\leq&
 \frac{C b_1\sigma^2}{(nb_0^d)^2}
 \sum_{i=1}^n
 \sum_{j=1\atop j \neq i}^{n}
 \frac{\mathds{1}
 (X_i \in\mathcal{X}_0)}
{\widehat{g}_{in}^2}
 K_0^2
\left(\frac{X_j-X_i}{b_0}\right)
\\
 &\leq&
 \frac{C b_1\sigma^2}{nb_0^d}
 \sum_{i=1}^n
 \frac{\mathds{1}
 (X_i\in \mathcal{X}_0)
 \widetilde{g}_{in}}
{\widehat{g}_{in}^2}\;,
 \label{VarSigmai}
\end{eqnarray}
where $\sigma^2=\Var(\varepsilon)$ and
$$
\widetilde{g}_{in}
 =
 \frac{1}{n b_0^d}
 \sum_{j=1,j \neq i}^{n}
 K_0^2\left( \frac{X_j-X_i}{b_0}\right).
$$
For the sum of  conditional covariances in (\ref{VarSigm}),
observe that by $(A_4)$ we have
\begin{eqnarray*}
 \lefteqn{
 \sum_{i=1}^n
 \sum_{j=1\atop j\neq i}^n
\Cov_n \left[
 \Sigma_{in}
 K_1^{(1)}
 \left(
 \frac{\varepsilon_{i}-\epsilon}{b_1}
 \right)
 ,
 \Sigma_{jn}
 K_1^{(1)}
 \left(
  \frac{\varepsilon_{j}-\epsilon}{b_1}
 \right)
\right]
 }
\\
&=&
 \sum_{i=1}^n
 \sum_{j=1\atop j\neq i}^n
\esp_n
 \left[
 \Sigma_{in}
 \Sigma_{jn}
 K_1^{(1)}
 \left(
\frac{\varepsilon_{i}-\epsilon}{b_1}
 \right)
 K_1^{(1)}
  \left(
\frac{\varepsilon_{j}-\epsilon}{b_1}
 \right)
  \right]
\\
& = &
 \sum_{i=1}^n
 \sum_{j=1\atop j\neq i}^n
 \frac{\mathds{1}(X_{i}\in\mathcal{X}_0)
\mathds{1}(X_{j}\in\mathcal{X}_0)}
 {(n b_0^d)^2 \widehat{g}_{in} \widehat{g}_{jn}}
 \sum_{k=1\atop k\neq i}^n
\sum_{\ell=1\atop \ell\neq j}^n
 K_0
 \left(
\frac{X_{k}-X_{i}}{b_0}
 \right)
 K_0
 \left(
\frac{X_{\ell}-X_{j}}{b_0} \right) \esp
 \left[
 \xi_{ki}
 \xi_{\ell j}
 \right],
\end{eqnarray*}
where
$$
\xi_{ki}
 =
 \varepsilon_k
 K_1^{(1)}
  \left(
\frac{\varepsilon_{i}-\epsilon}{b_1}
 \right).
$$
Moreover, under $(A_4)$, it is seen that for $k\neq\ell$,
$\esp[\xi_{ki}\xi_{\ell j}] =0 $ when $\Card\{i, j, k,
\ell\}\geq3$. Therefore
 the symmetry of $K_0$  yields that
\begin{eqnarray*}
 \lefteqn{
 \sum_{i=1}^n
 \sum_{j=1\atop j\neq i}^n
\Cov_n \left[ \Sigma_{in}
 K_1^{(1)}
  \left(
\frac{\varepsilon_{i}-\epsilon}{b_1}
 \right),
 \Sigma_{jn}
K_1^{(1)}
 \left(
 \frac{\varepsilon_{j}-\epsilon}{b_1}
 \right)
 \right]
 }
 &&
\\
& = & \sum_{i=1}^n
 \sum_{j=1\atop j\neq i}^n
\frac{\mathds{1}(X_{i} \in \mathcal{X}_0) \mathds{1} (X_{j}
 \in
\mathcal{X}_0)}
 {(n b_0^d)^2 \widehat{g}_{in} \widehat{g}_{jn}}
K_0^2 \left( \frac{X_{j}-X_{i}}{b_0} \right)
 \esp^2
 \left[
\varepsilon K_1^{(1)} \left( \frac{\varepsilon-\epsilon}{b_1}
\right) \right]
\\
&& + \sum_{i=1}^n
 \sum_{j=1\atop j\neq i}^n
 \frac{ \mathds{1} (X_{i}
\in \mathcal{X}_0) \mathds{1} (X_{j} \in \mathcal{X}_0)}
 {(nb_0^d)^2\widehat{g}_{in} \widehat{g}_{jn}}
 \sum_{k=1\atop k\neq  i, j}^n
 K_0 \left(\frac{X_{k}-X_{i}}{b_0}\right)
 K_0\left(\frac{X_{k}-X_{j}}{b_0}\right)
 \esp[\varepsilon^2]
\esp^2 \left[ K_1^{(1)}
 \left( \frac{\varepsilon-\epsilon}{b_1}
\right) \right].
\\
\end{eqnarray*}
Therefore, since
$$
\sup_{1\leq j\leq n}
 \left(
\frac{\mathds{1}\left(X_j\in\mathcal{X}_0\right)}
{|\widehat{g}_{jn}|} \right)
 =O_{\prob}(1)
$$
by Lemma \ref{Estig},  Lemma \ref{MomderK}-(\ref{MomderK1}) and
$(A_4)$ then give
\begin{eqnarray}
\nonumber
 \lefteqn{
  \left|
 \sum_{i=1}^n
 \sum_{j=1\atop j\neq i}^n
\Cov_n
 \left[
  \Sigma_{in}
 K_1^{(1)}
 \left(
\frac{\varepsilon_{i}-\epsilon}{b_1}
 \right)
 ,
  \Sigma_{jn}
 K_1^{(1)}
 \left(
 \frac{\varepsilon_{j}-\epsilon}{b_1}
\right) \right] \right|
 }
\\
  &=&
 O_{\prob}
 \left(
 \frac{b_1^4}{n b_0^d}
 \right)
 \sum_{i=1}^n
 \frac{\mathds{1}(X_i\in \mathcal{X}_0)\widetilde{g}_{in}}
 {|\widehat{g}_{in}|}
 +
 O_{\prob} (b_1^4)
 \sum_{i=1}^n
 \frac{\mathds{1}(X_i\in \mathcal{X}_0)|g_{in}|}
 {|\widehat{g}_{in}|}\;,
 \label{CovSigmaib}
\end{eqnarray}
where  $\widetilde{g}_{in}$ is defined as  in (\ref{VarSigmai})
and
$$
 g_{in}
 =
 \frac{1}{(n b_0^d)^2}
 \sum_{j=1\atop j\neq i}^n
 \sum_{k=1\atop k\neq j, i}^n
 K_0\left(\frac{X_k-X_i}{b_0}\right)
 K_0\left(\frac{X_k-X_j}{b_0}\right).
$$
The order of the first term in (\ref{CovSigmaib}) follows from
Lemma \ref{Estig}, which  gives
 \begin{eqnarray}
 \sum_{i=1}^n
 \frac{\mathds{1}(X_i\in \mathcal{X}_0)\widetilde{g}_{in}}
 {|\widehat{g}_{in}|}
 =
 O_{\prob}(n).
 \label{CovSigmaib1}
 \end{eqnarray}
Again, by Lemma \ref{Estig}, we have
\begin{eqnarray*}
 \sum_{i=1}^n
 \frac{\mathds{1}(X_i\in \mathcal{X}_0)|g_{in}|}
 {|\widehat{g}_{in}|}
 =
 O_{\prob}(1)
 \sum_{i=1}^n
 \mathds{1}
\left(X_i\in \mathcal{X}_0\right)|g_{in}|,
\end{eqnarray*}
with, using  the changes of variables $x_1=x_3+b_0z_1$,
$x_2=x_3+b_0z_2$,
\begin{eqnarray*}
\lefteqn{ \esp \left[ \sum_{i=1}^n \mathds{1}
\left(X_i\in\mathcal{X}_0\right)|g_{in}| \right] \leq
\frac{Cn^3}{(nb_0^d)^2} \esp \left|
K_0\left(\frac{X_3-X_1}{b_0}\right)
K_0\left(\frac{X_3-X_2}{b_0}\right) \right| }
\\
&\leq& \frac{Cn^3b_0^{2d}}{(nb_0^d)^2} \int\int\int
\left|K_0(z_1)K_0(z_2)\right| g(x_3+b_0z_1)g(x_3+b_0z_2)g(x_3)
dz_1dz_2dx_3.
\end{eqnarray*}
These bounds and the equality above, give under  $(A_2)$ and
$(A_7)$,
$$
 \sum_{i=1}^n
 \frac{\mathds{1}(X_i\in \mathcal{X}_0)|g_{in}|}
 {|\widehat{g}_{in}|}
 =
 O_{\prob}(n).
$$
Hence from  (\ref{CovSigmaib1}), (\ref{CovSigmaib}),
(\ref{VarSigmai}), (\ref{VarSigm}) and Lemma \ref{Estig}, we
deduce, for $b_1$ small enough,
\begin{eqnarray*}
\lefteqn{
 \Var_n
 \left[
 \sum_{i=1}^n
 \Sigma_{in}
 K_1^{(1)}
 \left(
 \frac{\varepsilon_i-\epsilon}{b_1}
 \right)
 \right]
 }
 &&
 \\
 & =&
 O_{\prob}
 \left(
 \frac{b_1}{n b_0^d}
 \right)
 \sum_{i=1}^{n}
 \frac{\mathds{1}(X_i\in \mathcal{X}_0)\widetilde{g}_{in}}
 {\widehat{g}_{in}^2}
 +
 O_{\prob}
 \left(
 \frac{b_1^4}{n b_0^d}
 \right)
\sum_{i=1}^n
 \frac{\mathds{1}(X_i \in \mathcal{X}_0)
\widetilde{g}_{in}} {|\widehat{g}_{in}|}
 +
 O_{\prob}(b_1^4)
 \sum_{i=1}^n
\frac{\mathds{1}(X_i \in\mathcal{X}_0)|g_{in}|}
{|\widehat{g}_{in}|}
\\
&=&
 O_{\prob}
\left( \frac{b_1}{b_0^d}
 +
 \frac{b_1^4}{b_0^d}
 +
 nb_1^4
 \right)
 =
 O_{\prob}
 \left(
 \frac{b_1}{b_0^d}
 +
 nb_1^4
 \right).
\end{eqnarray*}
Finally, this order, (\ref{EspSigmai}) and the Tchebychev
inequality give
$$
\sum_{i=1}^n
 \Sigma_{in}
 K_1^{(1)}
 \left(
 \frac{\varepsilon_i-\epsilon}{b_1}
 \right)
 =
 O_{\prob}
 \left(
 \frac{b_1}{b_0^d}+nb_1^4
 \right)^{1/2}. \eop
$$

\subsection*{Proof of Lemma \ref{BoundEspmchap}}

Define $\beta_{in}$  as in Lemma \ref{Betasum} and set
\begin{eqnarray*}
g_{in}
 =
 \frac{1}{nb_0^d}
 \sum_{j=1, j\neq i}^n
 K_0^4\left(\frac{X_j-X_i}{b_0}\right),
 \quad
\widetilde{g}_{in}
 =
 \frac{1}{nb_0^d}
 \sum_{j=1, j\neq i}^n
 K_0^2\left(\frac{X_j-X_i}{b_0}\right).
\end{eqnarray*}
The proof of the lemma  is based on the following bound:
\begin{eqnarray}
\esp_n
 \biggl[
 \mathds{1}
 \left(X_i\in\mathcal{X}_0\right)
 \left(\widehat{m}_{in} - m(X_i)\right)^k
\biggr]
 \leq
 C
 \left[
 \beta_{in}^k
 +
 \frac{
 \mathds{1}
 \left(X_i\in\mathcal{X}_0\right)
 \widetilde{g}_{in}^{k/2}}
 {(nb_0^d)^{(k/2)}\widehat{g}_{in}^k}
 \right],
 \quad
 k\in\{4,6\}.
 \label{Espm}
 \end{eqnarray}
 Indeed, taking successively $k=4$ and $k=6$ in (\ref{Espm}),
 we have, by (\ref{BetasumTBP}), Lemma \ref{Estig} and
 $(A_9)$,
\begin{eqnarray*}
 \sup_{1\leq i\leq n}
 \esp_n
 \biggl[
 \mathds{1}\left(X_i\in\mathcal{X}_0\right)
 \left(\widehat{m}_{in} - m(X_i)\right)^4
\biggr]
 &=&
 O_{\prob}
 \left(
 b_0^8
  +
 \frac{1}{(nb_0^d)^2}
 \right)
 =
 O_{\prob}
 \left(b_0^4+\frac{1}{nb_0^d}\right)^2,
 \\
 \sup_{1\leq i\leq n}
 \esp_n
 \biggl[
 \mathds{1}
 \left(X_i\in\mathcal{X}_0\right)
 \left(\widehat{m}_{in} - m(X_i)\right)^6
\biggr]
 &=&
 O_{\prob}
 \left(
 b_0^{12}
  +
 \frac{1}{(nb_0^d)^3}
 \right)
 =
 O_{\prob}
 \left(b_0^4+\frac{1}{nb_0^d}\right)^3,
\end{eqnarray*}
which gives the results of the Lemma.
 Hence it remains to prove (\ref{Espm}).
 For this,  define $\beta_{in}$ and  $\Sigma_{in}$  respectively as in
 Lemma \ref{Betasum} and Lemma \ref{Sigsum}.
 Since
 $\mathds{1}(X_i \in\mathcal{X}_0)
 \left(\widehat{m}_{in}- m(X_i)\right)
 =\beta_{in}+\Sigma_{in}$, and that
 $\beta_{in}$ depends only on
 $\left(X_1,\ldots,X_n\right)$, this gives, for $k\in\{4,6\}$
\begin{eqnarray}
\esp_n
 \biggl[
 \mathds{1}
 (X_i \in \mathcal{X}_0)
 \left(\widehat{m}_{in}- m(X_i)\right)^k
 \biggr]
  \leq
 C\beta_{in}^k
 +
 C\esp_n\left[\Sigma_{in}^k\right].
 \label{Espm5}
\end{eqnarray}
The order  of the second term of bound (\ref{Espm5}) is computed
by applying Theorem 2 in Whittle (1960) or the
Marcinkiewicz-Zygmund inequality (see e.g Chow and Teicher, 2003,
p. 386). These inequalities show that for linear form
$L=\sum_{j=1}^n a_j\zeta_j$ with independent mean-zero random
variables
 $\zeta_1,\ldots,\zeta_n$, it holds that, for any $k\geq 1$,
 $$
 \esp
 \left|L^k\right|
 \leq
  C(k)
  \left[
 \sum_{j=1}^n
  a_j^2
\esp^{2/k} \left|\zeta_j^k \right|
 \right]^{k/2},
 $$
where $C(k)$ is a positive real depending only on $k$. Now,
observe that for any $i\in[1,n]$,
$$
\Sigma_{in}
 =
 \sum_{j=1, j\neq i}^n
 \sigma_{jin},
 \quad
 \sigma_{jin}
 =
 \frac{\mathds{1}
 \left( X_i \in \mathcal{X}_0\right)}
 {nb_0^d\widehat{g}_{in}}
 \varepsilon_j
 K_0\left(\frac{X_j-X_i}{b_0}\right).
$$
Since under $(A_4)$, the $\sigma_{jin}$'s, $j\in[1,n]$, are
centered independent variables given $X_1,\ldots,X_n$,  this
yields, for any $k\in\{4,6\}$,
\begin{eqnarray*}
 \esp_n
 \left[
 \Sigma_{in}^k\right]
 \leq
 C\esp\left[\varepsilon^k\right]
 \left[
 \frac{\mathds{1}
 \left( X_i \in \mathcal{X}_0\right)}
 {(nb_0^d)^2\widehat{g}_{in}^2}
 \sum_{j=1}^n
  K_0^2
 \left(
 \frac{X_j-X_i}{b_0}
 \right)
 \right]^{k/2}
 \leq
 \frac{C\mathds{1}
 \left(X_i\in\mathcal{X}_0\right)\widetilde{g}_{in}^{k/2}}
 {(nb_0^d)^{(k/2)}\widehat{g}_{in}^k}\;.
\end{eqnarray*}
 Hence this  bound and
(\ref{Espm5})  give
$$
 \esp_n
 \biggl[
 \mathds{1}(X_i \in \mathcal{X}_0)
 \left(\widehat{m}_{in}- m(X_i)\right)^k
 \biggr]
 \leq C\left[
 \beta_{in}^k
  +
 \frac{
 \mathds{1}
 \left(X_i\in\mathcal{X}_0\right)
 \widetilde{g}_{in}^{k/2}}
 {(nb_0^d)^{(k/2)}\widehat{g}_{in}^k}
 \right],
$$
which proves (\ref{Espm}), and then completes the proof of the
lemma. \eop

\subsection*{Proof of Lemma \ref{Indep}}

Since $K_0(\cdot)$ has a compact support under $(A_7)$, there is a
$C>0$ such that $\| X_i - X_j \| \geq C b_0$ implies  that for any
integer number $k$ of $[1,n]$, $K_0 ( (X_k - X_i)/b_0) = 0$ if
$K_0 ( (X_j - X_k)/b_0) \neq 0$. Let $D_j \subset [1,n]$ be such
that an integer number $k$ of $[1,n]$ is in $D_j$ if and only if
$K_0 ( (X_j - X_k)/b_0) \neq 0$. Abbreviate $\prob (\cdot| X_1,
\ldots,X_n)$ into $\prob_n$ and assume that $\| X_i - X_j \| \geq
C b_0$ so that $D_i$ and $D_j$ have an empty intersection. Note
also that taking $C$ large enough ensures that $i$ is not in $D_j$
and $j$ is not in $D_i$. It then follows, under $(A_4)$ and since
$D_i$ and $D_j$ only depend upon $X_1,\ldots,X_n$,
\begin{eqnarray*}
\lefteqn{
 \prob_n
 \biggl(
 \left(
 \widehat{m}_{in} - m(X_i),\varepsilon_i
 \right)
 \in A \mbox{ \rm and }
 \left(
\widehat{m}_{jn} - m (X_j),\varepsilon_j \right)
 \in B \biggr)
 }
&&
\\
& = & \prob_n \left( \left( \frac{\sum_{k \in D_i \setminus\{i\}}
 \left( m(X_{k}) - m(X_i) + \varepsilon_{k} \right)
  K_0
\left( (X_{k} - X_i)/b_0 \right)} {\sum_{k \in D_i \setminus\{i\}}
 K_0 \left( (X_{k} - X_i)/b_0 \right)}
 ,
 \varepsilon_i\right) \in A \right.
\\
&& \;\;\;\;\;\;\;\;\;\;\;\;\;\;\;\;\;\;\;\; \left. \mbox{ \rm and}
\left( \frac{ \sum_{\ell \in D_j \setminus \{j \}}
 \left(
m(X_{\ell}) - m(X_j) + \varepsilon_{\ell} \right)
 K_0 \left(
(X_{\ell} - X_j)/b_0 \right) } { \sum_{\ell \in D_j \setminus \{j
\}} K_0 \left( (X_{\ell} - X_j)/b_0 \right)} ,
\varepsilon_j\right) \in B \right)
\\
& = & \prob_n \left( \left( \frac{ \sum_{k \in D_i \setminus \{i
\}} \left( m(X_{k}) - m(X_i) + \varepsilon_{k}\right)
 K_0 \left(
(X_{k} - X_i)/b_0 \right)} {\sum_{k \in D_i \setminus \{i \}}
K_0\left( (X_{k} - X_i)/b_0 \right)} ,
 \varepsilon_i \right) \in
A \right)
\\
&& \;\;\;\;\;\;\;\;\;\;\;\;\;\;\;\;\;\;\;\;
 \times\;
 \prob_n
 \left(
\left( \frac{ \sum_{\ell \in D_j \setminus \{j \}} \left(
m(X_{\ell}) - m(X_j) + \varepsilon_{\ell} \right) K_0 \left(
(X_{\ell} - X_j)/b_0 \right) } { \sum_{\ell \in D_j \setminus \{j
\}} K_0\left( (X_{\ell} - X_j)/b_0\right) } , \varepsilon_j
\right) \in B \right)
\\
& = & \prob_n
 \left(
 \left(\widehat{m}_{in} - m(X_i), \varepsilon_i \right)
 \in A
 \right)
\times \prob_n \left( \left(\widehat{m}_{jn} - m
(X_j),\varepsilon_j \right)
 \in B
\right).
\end{eqnarray*}
This gives the result of Lemma \ref{Indep}, since both
$\left(\widehat{m}_{in} - m (X_i), \varepsilon_i\right)$ and
$\left(\widehat{m}_{jn} - m (X_j), \varepsilon_j\right)$ are
independent given $X_1, \ldots, X_n$. \eop

\subsection*{Proof of Lemma \ref{sumzeta}}
Since $\widehat{m}_{in} - m(X_i)$ depends only upon
 $\left(X_1,\ldots,X_n,\varepsilon_k, k\neq i\right)$,
 we have
\begin{eqnarray*}
\sum_{i=1}^n
 \Var_n
 \left(\zeta_{in}\right)
  \leq
 \sum_{i=1}^n
 \esp_n
 \left[
 \zeta_{in}^2
 \right]
 =
 \sum_{i=1}^n
 \esp_n
 \left[
 \mathds{1}
 \left(X_i\in\mathcal{X}_0\right)
 \left(\widehat{m}_{in} - m(X_i)\right)^4
  \esp_{in}
  \left[
  K_1^{(2)}
  \left(\frac{\varepsilon_{i}-\epsilon}{b_1}\right)^2
  \right]
 \right],
 \end{eqnarray*}
with, using Lemma \ref{MomderK}-(\ref{MomderK2}),
\begin{eqnarray*}
 \esp_{in}
 \left[
 K_1^{(2)}
 \left(\frac{\varepsilon_{i}-\epsilon}{b_1}\right)^2
 \right]
 =
 \int
  K_1^{(2)}
 \left(
 \frac{e-\epsilon}{b_1}
 \right)^2
 f(e) de
 \leq
 C b_1.
\end{eqnarray*}
 Therefore these bounds and Lemma \ref{BoundEspmchap} give
\begin{eqnarray*}
\sum_{i=1}^n
 \Var_n\left(\zeta_{in}\right)
 &\leq&
 Cb_1
 \sum_{i=1}^n
 \esp_n
 \biggl[
 \mathds{1}
 \left(X_i\in\mathcal{X}_0\right)
 (\widehat{m}_{in} - m(X_i))^4
 \biggr]
 \\
 &\leq&
 Cnb_1
 \sup_{1\leq i\leq n}
 \esp_n
 \biggl[
 \mathds{1}
 \left(X_i\in\mathcal{X}_0\right)
 (\widehat{m}_{in} - m(X_i))^4
 \biggr]
 \\
 &\leq&
 O_{\prob}\left(nb_1\right)
 \left(b_0^4+\frac{1}{nb_0^d}\right)^2.
\end{eqnarray*}
which yields the desired result for the conditional variance.

\vskip 0.3cm
 We now prepare to compute the order of the conditional covariance.
 To that aim, observe that Lemma \ref{Indep} gives
\begin{eqnarray*}
\sum_{i=1}^n \sum_{j=1\atop j\neq i}^n \Cov_n \left(
\zeta_{in},\zeta_{jn}
 \right)
 =
  \sum_{i=1}^n
  \sum_{j=1\atop j\neq i}^n
  \mathds{1}
  \biggl(\left\|X_i - X_j\right\|<C b_0\biggr)
  \biggl(
  \esp_n
  \left[
  \zeta_{in}
  \zeta_{jn}
  \right]
  -
 \esp_n\left[\zeta_{in}\right]
 \esp_n\left[\zeta_{jn}\right]
 \biggr).
\end{eqnarray*}
 The order of the term above  is derived from the following
equalities:
\begin{eqnarray}
 \sum_{i=1}^n
 \sum_{j=1\atop j\neq i}^n
 \mathds{1}
 \biggl(
 \left\|X_i - X_j \right\|<C b_0
 \biggr)
 \esp_n\left[\zeta_{in}\right]
 \esp_n\left[\zeta_{jn}\right]
 &=&
 O_{\prob}
 \left(n^2b_0^db_1^6\right)
\left(b_0^4+ \frac{1}{nb_0^d}\right)^2,
  \label{Covzeta2}
  \\
  \sum_{i=1}^n
  \sum_{j=1\atop j\neq i}^n
  \mathds{1}
  \biggl(
  \left\|X_i - X_j\right\|<C b_0
  \biggr)
  \esp_n
  \left[
  \zeta_{in}
  \zeta_{jn}
  \right]
  &=&
  O_{\prob}
 \left(n^2b_0^db_1^{7/2}\right)
\left(b_0^4+ \frac{1}{nb_0^d}\right)^2.
 \label{Covzeta1}
\end{eqnarray}
 Indeed, since $b_1$ goes to $0$ under $(A_{10})$,
 (\ref{Covzeta2}) and (\ref{Covzeta1}) yield that
\begin{eqnarray*}
 \sum_{i=1}^n
 \sum_{j=1\atop j\neq i}^n
 \Cov_n
 \left(\zeta_{in},\zeta_{jn}\right)
 &=&
O_{\prob}
 \left[
 \left(n^2b_0^db_1^6\right)
 \left(b_0^4+ \frac{1}{nb_0^d}\right)^2
 +
 \left(n^2b_0^db_1^{7/2}\right)
 \left(b_0^4+ \frac{1}{nb_0^d}\right)^2
\right]
\\
&=& O_{\prob}
 \left(n^2b_0^db_1^{7/2}\right)
 \left(b_0^4+ \frac{1}{nb_0^d}\right)^2,
\end{eqnarray*}
 which gives the result for the conditional
covariance. Hence,  it remains to prove (\ref{Covzeta2}) and
(\ref{Covzeta1}).  For (\ref{Covzeta2}), note that by $(A_4)$ and
Lemma \ref{MomderK}-(\ref{MomderK2}), we have
 \begin{eqnarray*}
 \left|
 \esp_{n}
 \left[\zeta_{in}\right]
 \right|
 &=&
 \left|
 \esp_n
 \left[
 \mathds{1}
 \left(X_i\in\mathcal{X}_0\right)
 (\widehat{m}_{in} - m(X_i))^2
 \esp_{in}
 \left[
 K_1^{(2)}
 \left(
 \frac{\varepsilon_i-\epsilon}{b_1}
 \right)
 \right]
 \right]
 \right|
 \\
 &\leq&
 Cb_1^3
 \biggl(
  \esp_n
 \biggl[
 \mathds{1}
 \left(X_i\in\mathcal{X}_0\right)
 (\widehat{m}_{in} - m(X_i))^4
 \biggr]
 \biggr)^{1/2}.
 \end{eqnarray*}
Hence from this bound and Lemma \ref{BoundEspmchap} we deduce
\begin{eqnarray*}
\sup_{1\leq i, j\leq n}
 \left|
 \esp_{n}
 \left[\zeta_{in}\right]
 \esp_n
\left[\zeta_{jn}\right]
 \right|
 &\leq&
 Cb_1^6
 \sup_{1\leq i \leq n}
 \esp_n
 \biggl[
 \mathds{1}
 \left(X_i\in\mathcal{X}_0\right)
 (\widehat{m}_{in} - m(X_i))^4
 \biggr]
 \\
 &\leq&
 O_{\prob}
 \left(b_1^6\right)
\left(b_0^4+ \frac{1}{nb_0^d}\right)^2.
\end{eqnarray*}
Therefore, since the Markov inequality gives
\begin{eqnarray}
\sum_{i=1}^n
 \sum_{j=1\atop j\neq i}^n
 \mathds{1}
 \biggl(
 \| X_i - X_j \|< C b_0
 \biggr)
 =
 O_{\prob}(n^2 b_0^d),
 \label{Markov}
 \end{eqnarray}
 it then follows that
\begin{eqnarray*}
  \sum_{i=1}^n
 \sum_{j=1\atop j\neq i}^n
 \mathds{1}
 \biggl(
 \| X_i - X_j \|< C b_0
 \biggr)
 \esp_{n}
 \left[\zeta_{in}\right]
 \esp_{n}
 \left[\zeta_{jn}\right]
  =
 O_{\prob}
 \left(n^2b_0^db_1^6\right)
\left(b_0^4+ \frac{1}{nb_0^d}\right)^2,
\end{eqnarray*}
which   proves (\ref{Covzeta2}).

\vskip 0.3cm
 For (\ref{Covzeta1}), set $Z_{in}= \mathds{1}
\left(X_i\in\mathcal{X}_0\right)\left(\widehat{m}_{in} -
m(X_i)\right)^2$, and note that for $i\neq j$, we have
 \begin{eqnarray}
 \esp_n\left[\zeta_{in}\zeta_{jn}\right]
 =
 \esp_n
 \left[
 Z_{in}
  K_1^{(2)}
 \left(
 \frac{\varepsilon_j-\epsilon}{b_1}
 \right)
  \esp_{in}
 \left[
 Z_{jn}
  K_1^{(2)}
 \left(
 \frac{\varepsilon_i-\epsilon}{b_1}
 \right)
 \right]
 \right],
 \label{Prodzeta}
 \end{eqnarray}
where
\begin{eqnarray}
\nonumber
 \lefteqn{
 \esp_{in}
 \left[
 Z_{jn}
  K_1^{(2)}
 \left(
 \frac{\varepsilon_i-\epsilon}{b_1}
 \right)
 \right]
}
\\\nonumber
&=&
 \beta_{jn}^2
 \esp_{in}
 \left[
 K_1^{(2)}
 \left(
 \frac{\varepsilon_i-\epsilon}{b_1}
 \right)
 \right]
 +
 2\beta_{jn}
 \esp_{in}
 \left[
 \Sigma_{jn}
 K_1^{(2)}
 \left(
 \frac{\varepsilon_i-\epsilon}{b_1}
 \right)
 \right]
 +
 \esp_{in}
 \left[
 \Sigma_{jn}^2
 K_1^{(2)}
 \left(
 \frac{\varepsilon_i-\epsilon}{b_1}
 \right)
 \right].
 \\
\label{Covzeta3}
\end{eqnarray}
The first term of Equality (\ref{Covzeta3}) is treated by using
Lemma \ref{MomderK}-(\ref{MomderK2}). This gives
\begin{eqnarray}
\left| \beta_{jn}^2
 \esp_{in}
 \left[
 K_1^{(2)}
 \left(
 \frac{\varepsilon_i-\epsilon}{b_1}
 \right)
 \right]
 \right|
 \leq
 Cb_1^3
 \beta_{jn}^2.
 \label{Covzeta4}
\end{eqnarray}
Since  under $(A_4)$, the $\varepsilon_j$'s are independent
centered variables, and  are independent of the  $X_j$'s, the
second term in (\ref{Covzeta3}) gives
\begin{eqnarray*}
\esp_{in}
 \left[
 \Sigma_{jn}
   K_1^{(2)}
 \left(
 \frac{\varepsilon_i-\epsilon}{b_1}
 \right)
 \right]
&=&
 \frac{\mathds{1}\left(X_j\in\mathcal{X}_0\right)}
 {nb_0^d\widehat{g}_{jn}}
\sum_{k=1, k\neq j}^n
 K_0
 \left(
 \frac{X_k-X_j}{b_0}
 \right)
\esp_{in}
 \left[
 \varepsilon_k
   K_1^{(2)}
 \left(
 \frac{\varepsilon_i-\epsilon}{b_1}
 \right)
 \right]
 \\
 &=&
 \frac{\mathds{1}\left(X_j\in\mathcal{X}_0\right)}
 {nb_0^d\widehat{g}_{jn}}
 K_0
 \left(
 \frac{X_i-X_j}{b_0}
 \right)
\esp_{in}
 \left[
 \varepsilon_i
   K_1^{(2)}
 \left(
 \frac{\varepsilon_i-\epsilon}{b_1}
 \right)
 \right].
\end{eqnarray*}
Therefore, by $(A_7)$ which ensures that $K_0$ is bounded, the
equality above and Lemma \ref{MomderK}-(\ref{MomderK2}) yield that
\begin{eqnarray}
\left| \beta_{jn} \esp_{in}
 \left[
 \Sigma_{jn}
   K_1^{(2)}
 \left(
 \frac{\varepsilon_i-\epsilon}{b_1}
 \right)
 \right]
\right|
 \leq
 Cb_1^3
\left| \beta_{jn}
\frac{\mathds{1}\left(X_j\in\mathcal{X}_0\right)}
 {nb_0^d\widehat{g}_{jn}}
\right|.
 \label{Covzeta5}
\end{eqnarray}
For the  last term in (\ref{Covzeta3}), we have
\begin{eqnarray*}
\lefteqn{
 \esp_{in}
 \left[
 \Sigma_{jn}^2(x)
   K_1^{(2)}
 \left(\frac{\varepsilon_i-\epsilon}{b_1}\right)
 \right]
}
\\
&=& \frac{1}{(nb_0^d\widehat{g}_{jn})^2}
 \sum_{k=1\atop k\neq j}^n
\sum_{\ell=1\atop\ell\neq j}^n
 K_0\left(\frac{X_k-X_j}{b_0}\right)
 K_0\left(\frac{X_{\ell}-X_j}{b_0}\right)
 \esp_{in}
 \left[
 \varepsilon_k
 \varepsilon_{\ell}
  K_1^{(2)}
 \left(\frac{\varepsilon_i-\epsilon}{b_1}\right)
 \right]
 \\
 &=&
 \frac{1}{(nb_0^d\widehat{g}_{jn})^2}
 \sum_{k=1, k\neq j}^n
 K_0^2\left(\frac{X_k-X_j}{b_0}\right)
 \esp_{in}
 \left[
 \varepsilon_k^2
  K_1^{(2)}
 \left(\frac{\varepsilon_i-\epsilon}{b_1}\right)
 \right],
\end{eqnarray*}
with, using Lemma \ref{MomderK}-(\ref{MomderK2}),
\begin{eqnarray*}
\lefteqn
 {
 \left|
 \esp_{in}
 \left[
 \varepsilon_k^2
  K_1^{(2)}
 \left(\frac{\varepsilon_i-\epsilon}{b_1}\right)
 \right]
 \right|
}
\\
&\leq& \max
 \left\lbrace
\sup_{e\in\Rit} \left|
 \esp_{in}
 \left[
 \varepsilon^2
  K_1^{(2)}
 \left(\frac{\varepsilon-e}{b_1}\right)
 \right]
 \right|,
 \;
 \esp[\varepsilon^2]
\sup_{e\in\Rit}
 \left|
 \esp_{in}
 \left[
 K_1^{(2)}
 \left(\frac{\varepsilon-e}{b_1}\right)
 \right]
 \right|
\right\rbrace
\\
&\leq&
 C b_1^3.
\end{eqnarray*}
Therefore
$$
 \left|
 \esp_{in}
 \left[
 \Sigma_{jn}^2
   K_1^{(2)}
 \left(
 \frac{\varepsilon_i-\epsilon}{b_1}
 \right)
 \right]
 \right|
 \leq
\frac{Cb_1^3}{(nb_0^d\widehat{g}_{jn})^2} \sum_{k=1, k\neq j}^n
K_0^2\left(\frac{X_k-X_j}{b_0}\right).
$$
Substituting this bound, (\ref{Covzeta5}) and (\ref{Covzeta4}) in
(\ref{Covzeta3}), we obtain
$$
\left|
 \esp_{in}
 \left[
 Z_{jn}
  K_1^{(2)}
 \left(
 \frac{\varepsilon_i-\epsilon}{b_1}
 \right)
 \right]
\right|
 \leq
 Cb_1^3 M_n,
$$
where
$$
M_n = \sup_{1\leq j\leq n} \left[ \beta_{jn}^2 + \left| \beta_{jn}
\frac{\mathds{1}\left(X_j\in\mathcal{X}_0\right)}
{nb_0^d\widehat{g}_{jn}} \right| +
\frac{1}{(nb_0^d\widehat{g}_{jn})^2} \sum_{k=1, k\neq j}^n
K_0^2\left(\frac{X_k-X_j}{b_0}\right) \right].
$$
Hence from (\ref{Prodzeta}), the Cauchy-Schwarz inequality, Lemma
\ref{BoundEspmchap} and Lemma \ref{MomderK}-(\ref{MomderK2}), we
deduce
\begin{eqnarray*}
\lefteqn{ \sum_{i=1}^n \sum_{j=1\atop j\neq i}^n
 \mathds{1}
 \biggl(
 \| X_i - X_j \|< C b_0
 \biggr)
\left| \esp_n \left[ \zeta_{in}\zeta_{jn} \right] \right| }
\\
&\leq& C M_nb_1^3 \sum_{i=1}^n \sum_{j=1\atop j\neq i}^n
 \mathds{1}
 \biggl(
 \| X_i - X_j \|< C b_0
 \biggr)
\esp_n
 \left|
 Z_{in} K_1^{(2)}
\left(\frac{\varepsilon_j-\epsilon}{b_1}\right) \right|
\\
&\leq& CM_nb_1^3 \sum_{i=1}^n \sum_{j=1\atop j\neq i}^n
 \mathds{1}
 \biggl(
 \| X_i - X_j \|< C b_0
 \biggr)
\esp_n^{1/2} \left[ Z_{in}^2 \right]
 \esp_n^{1/2}
 \left[
 K_1^{(2)}
 \left(
 \frac{\varepsilon_j-\epsilon}{b_1}
 \right)^2
 \right]
 \\
 &\leq&
 M_n b_1^3
 O_{\prob}
 \left(b_0^4+\frac{1}{nb_0^d}\right)
 (b_1)^{1/2}
\sum_{i=1}^n \sum_{j=1\atop j\neq i}^n \biggl( \mathds{1}
\left(\|X_i-X_j\|\leq Cb_0\right) \biggr).
\end{eqnarray*}
Moreover, (\ref{BetasumTBP}) and Lemma \ref{Estig} give, under
$(A_1)$, $(A_7)$ and $(A_9)$,
$$
M_n =
 O_{\prob}
 \left( b_0^4 + \frac{b_0^2}{nb_0^d}
  +
\frac{1}{nb_0^d}
 \right)
= O_{\prob} \left(b_0^4+ \frac{1}{nb_0^d}\right).
$$
Finally, substituting this order in the bound above, and using
(\ref{Markov}), we arrive at
\begin{eqnarray*}
\sum_{i=1}^n \sum_{j=1\atop j\neq i}^n
 \mathds{1}
 \biggl(
 \| X_i - X_j \|< C b_0
 \biggr)
\esp_n \left[\zeta_{in}\zeta_{jn}\right] =
 O_{\prob}
\left(n^2b_0^db_1^{7/2}\right)
\left(b_0^4+\frac{1}{nb_0^d}\right)^2.
\end{eqnarray*}
This proves (\ref{Covzeta1}), and then completes the proof of the
theorem. \eop

\end{document}